\documentclass[aip,superscriptaddress]{revtex4}
\usepackage{amsmath,amssymb,amsfonts}
\usepackage{bbold}
\usepackage{color}
\usepackage{url}
\usepackage{graphicx}
\usepackage{tikz}
\usepackage{epstopdf}

\usepackage{hyperref}

\newcommand{\diag}{{\rm diag\,}}

\newcommand{\tr}{{\rm tr\,}}

\newcommand{\eins}{\leavevmode\hbox{\small1\kern-3.8pt\normalsize1}}

\begin{document}

\newtheorem{punkt}{}[section]

\newtheorem{corollary}[punkt]{Corollary}
\newtheorem{lemma}[punkt]{Lemma}
\newtheorem{proposition}[punkt]{Proposition}
\newtheorem{theorem}[punkt]{Theorem}
\newtheorem{remark}[punkt]{Remark}
\newtheorem{remarks}[punkt]{Remarks}
\newtheorem{example}[punkt]{Example}
\newtheorem{examples}[punkt]{Examples}
\newtheorem{problem}[punkt]{Problem}
\newtheorem{problems}[punkt]{Problems}
\newtheorem{question}[punkt]{Question}
\newtheorem{questions}[punkt]{Questions}
\newtheorem{definition}[punkt]{Definition}
\newtheorem{conjecture}[punkt]{Conjecture}
\newtheorem{assumption}[punkt]{Assumption}
\newtheorem{assumptions}[punkt]{Assumptions}
\newtheorem{construction}[punkt]{Construction}

\title[Additive Matrix Convolutions]{Additive Matrix Convolutions of P\'olya Ensembles and Polynomial Ensembles}
\author{Mario Kieburg}\email[]{mkieburg@physik.uni-bielefeld.de}
\affiliation{Fakult\"at f\"ur Physik,
  Universit\"at Bielefeld, Bielefeld, Germany}

\newcommand{\corr}[1]{{\color{red}#1}}
\newcommand{\tim}[1]{{\color{blue}#1}}


\begin{abstract}
Recently subclasses of polynomial ensembles for additive and multiplicative matrix convolutions were identified which were called P\'olya ensembles (or polynomial ensembles of derivative type). Those ensembles are closed under the respective convolutions and, thus, build a semi-group when adding by hand a unit element. They even have a semi-group action on the polynomial ensembles.  Moreover in several works transformations of the bi-orthogonal functions and kernels of a given polynomial ensemble  were derived when performing an additive or multiplicative matrix convolution with particular P\'olya ensembles. For the multiplicative matrix convolution on the complex square matrices the transformations were even done for general P\'olya ensembles. In the present work we generalize these results to the additive convolution on Hermitian matrices, on Hermitian anti-symmetric matrices, on Hermitian anti-self-dual matrices and on rectangular complex matrices. For this purpose we derive the bi-orthogonal functions and the corresponding kernel for a general P\'olya ensemble which was not done before. With the help of these results we find transformation formulas for the convolution with a fixed matrix or a random matrix drawn from a general polynomial ensemble. As an  example we consider P\'olya ensembles with an associated weight which is a P\'olya frequency function of infinite order. But we also explicitly evaluate the Gaussian unitary ensemble as well as the complex Laguerre (aka Wishart, Ginibre or chiral Gaussian unitary) ensemble. All results hold for finite matrix dimension. Furthermore we derive a recursive relation between Toeplitz determinants which appears as a by-product of our results.
\end{abstract}

\maketitle

\section{Introduction}\label{sec:intro}

Convolutions on matrix spaces have a long tradition in mathematics, physics and beyond. Very early products and sums of random matrices were studied either because of group theoretical interest, see~\cite{Helgason:2000} for a textbook on harmonic analysis on Lie groups and references therein, to generalize the central limit theorem to non-commutative operators~\cite{Bellman:1954}, to generlize random walks on operators~\cite{Dyson:1962} or to study stability problems~\cite{May:1972}. These models can be interpreted as discrete stochastic processes and found a rich variety of applications, to name only a few: ecological systems (see \cite{Allesina:2015} for a recent review), condensed matter physics (see~\cite{Beenakker:1997} for a review), and wireless telecommunication (see~\cite{Tulino:2004}). See also~\cite{Schehr:2017} in which several chapters deal with stochastic processes on matrix spaces. Dyson's Brownian motion~\cite{Dyson:1962} is a prominent example where matrix convolutions play an important role. There is a vast literature, e.g. see~\cite[Chapter 3]{Katori:2016} and references therein, which only deals with the spectral statistics of Dyson's Brownian motion and its applications. The problem in all these models of products and sums of random matrices was that they were very restrictive to particular probability distributions of random matrices like the Gaussian to derive all spectral properties at finite matrix dimension.

The situation changed in the past years. The interest in the spectral statistics of products of finite dimensional random matrices triggered a revival of these old problems, see~\cite{Akemann:2012} for the first works on the complex eigenvalues and~\cite{Akemann:2013} for the first work on the singular values on this new development. The new approaches led to the development of new techniques, e.g. see~\cite{Akemann:2015} for a reviews and~\cite{Kieburg:2016a,Kieburg:2016b} for new conceptual applications of harmonic analysis, and nurtured ideas to apply those tools also to sums of random matrices~\cite{Kuijlaars:2016a}. One of these  technical tools is the concept of polynomial ensembles~\cite{Kuijlaars:2014}, see Definition~\ref{def:PolEns}. A random matrix $X$ drawn from a polynomial ensemble $P(X)$ has a joint probability density $p(a)$ of its eigenvalues or its singular values $a$, depending on the considered matrix $X$, which has the form
\begin{equation}\label{jpdf-pol-ens.intro}
	p(a)\propto\Delta_n(a)\det[w_b(a_c)]_{b,c=1,\ldots,n}
\end{equation}
with  $\Delta_n(a)$ the Vandermonde determinant and $w_b$ some functions. There are two advantages of such ensembles. One is that they correspond to determinantal point processes, being a particular form of a bi-orthogonal ensemble~\cite{Borodin:1999}. The second advantage becomes immediate when one asks for the spectral statistics of a product $XY$ or a sum $X+Y$ with another statistically independent random matrix $Y$ drawn from the density $\tilde{P}(Y)$. As in the univariate case one has to perform a convolution of the probability densities $P$ and $\tilde{P}$ either of multiplicative or additive kind. In some cases the group integrals involved in such a convolution can be evaluated. These group integrals are of the form like the Harish-Chandra-Itzykson-Zuber integral~\cite{Harish:1956,Itzykson:1980} or the Berezin-Karpelevic integral~\cite{Berezin:1958,Guhr:1996} and are ratios of determinants with the Vandermonde determinant $\Delta_n(a)$ in the denominator. This is the point where the particular form~\eqref{jpdf-pol-ens.intro} of a polynomial ensemble comes into the game. Only due to the Vandermonde determinant in Eq.~\eqref{jpdf-pol-ens.intro} a further computation is possible since it cancels with the one from the group integral, see \cite{Akemann:2013,Akemann:2015,Kieburg:2016a,Kieburg:2016b,Kuijlaars:2014,Kuijlaars:2016a,Kieburg:2017} for several examples.

The theoretical development of matrix convolutions did not stop with polynomial ensembles. The reason is that the sum or the product of two polynomial ensembles does not necessarily yield a polynomial ensemble. Hence the form~\eqref{jpdf-pol-ens.intro} would be immediately lost after one ``time" step in an additive or multiplicative stochastic process. Recently subclasses of polynomial ensembles were identified which were closed under these matrix convolutions, see \cite{Kieburg:2016a,Kieburg:2016b} for the multiplication of complex square matrices, \cite{Kuijlaars:2016a} for the summation of Hermitian matrices, \cite{Kieburg:2017} for the summation of Hermitian antisymmetric matrices, Hermitian anti-self-dual matrices and complex rectangular matrices. All these matrix spaces where these subclasses were found belong to the Dyson index $\beta=2$ \cite{Zirnbauer:1996,Altland:1997}. These subclasses were first named polynomial ensembles of derivative type~\cite{Kieburg:2016a,Kieburg:2016b,Kuijlaars:2016a} since the functions have the form $w_b(x)=D^{b-1}\omega(x)$ with $D$ a differential operator.  In~\cite{Kieburg:2017} the name P\'olya ensemble was proposed since the corresponding weights are related to P\'olya frequency functions, see below. The differential operator $D$ depends on the space of matrices and the kind of convolution. Interestingly $D$ is a differential operator of first order for the multiplicative convolution on complex square matrices~\cite{Kieburg:2016a,Kieburg:2016b} and the additive convolution of Hermitian matrices~\cite{Kuijlaars:2016a}. However $D$ is of second order for the additive convolution on Hermitian matrices, on Hermitian anti-symmetric matrices, on Hermitian anti-self-dual matrices and on rectangular complex matrices and is structural of the same form only depending on an index $\nu=\{\pm1/2\}\cup\mathbb{N}_0$, see~\cite{Kieburg:2017} and Definition~\ref{def:PolEns}. The index $\nu$ relates to the level repulsion from the origin and is well-known in the case of Laguerre ensembles~\cite{Mehta:2004,Forrester:2010,Akemann:2011}, especially that these three matrix spaces can be dealt in a unifying way.

A question remained to be answered, namely which functions $\omega$ can be chosen in a P\'olya ensemble such that~\eqref{jpdf-pol-ens.intro} is still a probability density. In a very recent work~\cite{Kieburg:2017} this question was addressed and the suitable functions were related to P\'olya frequency functions~\cite{Polya:1913,Polya:1915,Schoenberg:1951,Karlin:1968}. A P\'olya frequency functions $f$ of order $n\in\mathbb{N}$ on $\mathbb{R}$ satisfies the positivity conditions~\cite{Polya:1913,Polya:1915,Schoenberg:1951,Karlin:1968}
\begin{equation}\label{def:Polyafreq}
\Delta_j(x)\Delta_j(y)\det[f(x_b-y_c)]_{b,c=1,\ldots,j}\geq0,\ {\rm for\ all}\ x,y\in\mathbb{R}^j\ {\rm and}\ j=1,\ldots,n.
\end{equation}
When this inequality holds for all $n\in\mathbb{N}$ the function $f$ is called a P\'olya frequency function of infinite order and has a particularly simple and explicit form in terms of its Fourier transform~\cite{Schoenberg:1951,Karlin:1968}, see also Example~\ref{example.a}.

After the suitable subclasses of polynomial ensembles were identified the question of their general spectral statistics comes into mind. Due to the particular form~\eqref{jpdf-pol-ens.intro} one usually thinks of the bi-orthogonal functions and the kernels of the corresponding determinantal point processes, see~\cite{Borodin:1999} for the general approach with bi-orthogonal ensembles. For general P\'olya ensembles corresponding to the multiplicative convolution on complex square matrices this question was recently answered in~\cite{Kieburg:2016b}. In the present work we will generalize these results to general P\'olya ensembles corresponding to the additive convolution on Hermitian matrices, on Hermitian anti-symmetric matrices, on Hermitian anti-self-dual matrices and on rectangular complex matrices  (or in short on $H_2={\rm Herm}(n)$, $H_1=\imath{\rm o}(n)$, $H_4=\imath{\rm usp}(2n)$ and $M_\nu=\mathbb{C}^{n\times(n+\nu)}$, respectively). This will be our first main result. For the Gaussian unitary ensembles (GUE) and the complex Laguerre ensemble, these results readily reduce to the known results~\cite{Mehta:2004,Forrester:2010,Akemann:2011}.

Another question which will be addressed in the present work is regarding the statistics when a general P\'olya ensembles on $H_2$, $H_1$, $H_4$ and $M_\nu$ is shifted by a constant matrix in the same space. We again derive explicit expression for their bi-orthogonal functions and kernels. In particular cases, like for Gaussian probability densities, this is already known~\cite{Claeys:2015}.

The third result will be the generalization of the transformation formulas when we add to a general polynomial ensembles in one of the spaces $H_2$, $H_1$, $H_4$ and $M_\nu$ with given bi-orthogonal functions and kernel a P\'olya ensemble. For the multiplicative convolution this was already done for the Gaussian case and the case of truncated unitary matrices (Jacobi ensemble) in~\cite{Kieburg:2015,Kuijlaars:2015,Kuijlaars:2016a} and for general P\'olya ensembles in~\cite{Kieburg:2016b}. A similar approach was employed in~\cite{Forrester:2017} where the natural action of the general linear group ${\rm Gl}_{\mathbb{C}}(n)$ distributed by the induced Ginibre ensemble on the Hermitian matrices distributed by a general polynomial ensemble was considered. For the additive convolution with the GUE and the complex Laguerre ensemble such transformation fomulas were recently derived in~\cite{Claeys:2015,Kuijlaars:2016a}. The authors of~\cite{Kuijlaars:2016a} related the problem to Gelfand pairs, $(G,K)$ with $K$ a compact subgroup of the Lie group $G$. For the additive convolution on Hermitian matrices the Gelfand pairs is $(G,K)=({\rm U}(n)\ltimes{\rm Herm}(n),{\rm U}(n)\ltimes\{0\})$ with the semi-direct product on $G={\rm U}(n)\ltimes{\rm Herm}(n)$ given by~\cite[Sec.~2.1]{Kuijlaars:2016a}
\begin{equation}\label{semi-direct}
(U_1,H_1)\cdot(U_2,H_2)=(U_1U_2,H_1+U_1H_2U_1^*)
\end{equation}
$U_2^*$ the Hermitian adjoint of $U_2$. In this framework we consider in the present work the four kinds of Gelfand pairs $({\rm O}(n)\ltimes\imath{\rm o}(n),{\rm O}(n)\ltimes\{0\})$, $({\rm U}(n)\ltimes{\rm Herm}(n),{\rm U}(n)\ltimes\{0\})$, $({\rm USp}(2n)\ltimes\imath{\rm usp}(2n),{\rm USp}(2n)\ltimes\{0\})$ and $(({\rm U}(n)\times{\rm U}(n+\nu))\ltimes\mathbb{C}^{n\times(n+\nu)},({\rm U}(n)\times{\rm U}(n+\nu))\ltimes\{0\})$ with a similar semi-direct product as Eq.~\eqref{semi-direct}. 

The work is built up as follows. In Sec.~\ref{sec:pre} we introduce our notation and the setting which we consider. In particular we state the starting points of our study which are the joint probability densities of P\'olya ensembles on $H_2$, $H_1$, $H_4$ and $M_\nu$ without a shift (Lemma~\ref{lem:jpdf.Polya}), with a shift with a constant matrix (Theorem~\ref{thm:jpdf.Polya.fixed}), and with a shift with a general polynomial ensemble (Theorem~\ref{thm:jpdf.Polya.poly}). In Sec.~\ref{sec:H2} we only derive the bi-orthogonal functions and kernels corresponding to the additive convolution on the Hermitian matrices, $H_2$. As a by-product we derive a recursive relation between Toeplitz determinants, see Corollary~\ref{cor:toeplitz}.  The case of the additive convolution on $H_1$, $H_4$ and $M_\nu$ is considered in Sec.~\ref{sec:M}. As already pointed out the three matrix spaces $H_1$, $H_4$ and $M_\nu$ can be dealt in a unifying way with a parameter $\nu$ which is $\nu=\pm1/2$ for $H_1$, $\nu=+1/2$ for $H_4$, and $\nu\in\mathbb{N}_0$ for $M_\nu$. In Sec.~\ref{sec:conclusio} we summarize our results. All results are given for finite matrix dimension. The asymptotic analysis of the P\'olya ensembles is not the aim of the present work.

\section{Preliminaries}\label{sec:pre}

We consider the additive convolution on either one of the three classical compact Lie algebras times the imaginary unit (Hermitian antisymmetric ($\beta=1$, $H_{1}=\imath\,{\rm o}(n)$), Hermitian ($\beta=2$, $H_{2}={\rm Herm}(n)$) or Hermitian anti-self-dual ($\beta=4$, $H_{4}=\imath\,{\rm usp}(2n)$) matrices), and of complex $n\times (n+\nu)$ matrices which can be cast into the chiral form
\begin{equation}
M_{\nu}=\left\{\left.\left[\begin{array}{cc} 0 & W \\ W^* & 0 \end{array}\right]\right|W\in{\rm Mat}_{\mathbb{C}}(n,n+\nu)\right\}.
\end{equation}
Here, we employ the notation of the work~\cite{Kieburg:2017}. The corresponding compact groups keeping these spaces invariant under their adjoint action are the three classical groups (orthogonal ($\beta=1$, $K_1={\rm O}(n)$), unitary ($\beta=2$, $K_2={\rm U}(n)$) or unitary symplectic ($\beta=4$, $K_4={\rm USp}(2n)$) matrices), and the group $\hat{K}_{\nu}={\rm U}(n)\times{\rm U}(n+\nu)$. The indices $\beta$ and $\nu$ are also known as the Dyson index and the topological charge. $L^1$-functions on one of the sets $M=H_\beta,M_\nu$ are called $\mathcal{K}$-invariant with $\mathcal{K}=K_\beta,\hat{K}_\nu$, respectively, are defined as
\begin{equation}\label{group-inv-func}
L^{1,\mathcal{K}}(M)=\left\{\left.f_{M}\in L^{1}(M)\right| f_{M}(km k^*)=f_{M}(m)\ \forall m\in M,\,k\in\mathcal{K}\right\}.
\end{equation}
We denote the Hermitian adjoint of a matrix $k$ by $k^*$.

$\mathcal{K}$-invariant functions only depend on the eigenvalues for $M=H_2$ or on their squared singular values for $M=H_1,H_4,M_\nu$. Thus we need the space of diagonal real $n\times n$ matrices $D\simeq\mathbb{R}^n$ and of diagonal positive definite $n\times n$ matrices: $A=\exp[D]\simeq\mathbb{R}_+^n$. Furthermore the $\mathcal{K}$-invariance of a function $f_M\in L^{1,\mathcal{K}}(M)$ carries over to an invariance under the symmetric group $\mathbb{S}$ of $n$ elements for the corresponding function $f_\mathcal{D}\in L^{1,\mathbb{S}}(\mathcal{D})$ with $\mathcal{D}=D,A$, respectively.

We equip the matrix spaces $H_\beta,M_\nu,D$ and $A$ with the flat Lebesgue measures denoted by $dy$, $dg$, $da$ etc. and the groups $K_\beta$ and $\hat{K}_\nu$ with the normalized Haar measure denoted by $d^*k$. The relations between $f_M\in L^{1,\mathcal{K}}(M)$ and $f_\mathcal{D}\in L^{1,\mathbb{S}}(\mathcal{D})$ are given by the isometries (with respect to the $L^1$-norm \mbox{$\|\,\cdot\,\|_1$}):
 \begin{equation}\label{I-H2}
 \begin{split}
 \mathcal{I}_{H_2}:&\, L^{1,K_2}(H_2)\rightarrow L^{1,{\mathbb{S}}}(D),\quad f_D(a)=\mathcal{I}_{H_2} f_{H_2}(a)=C_n f_{H_2}(a) \Delta_n^2(a), \quad a \in D,\\
 \mathcal{I}_{M} :&\, L^{1,\mathcal{K}}(M) \to L^{1,{\mathbb{S}}}(A),\quad f_A(a)=\mathcal{I}_{M}f_{M}(a)=C_{n,\nu}^*\det a^\nu f_{M_\nu}\left(\iota_M(a)\right)\Delta_n^2(a),\quad a\in A,
 \end{split}
 \end{equation}
 where $\nu\in\mathbb{N}_0$ for $(M,\mathcal{K})=(M_\nu,\hat{K}_\nu)$, $\nu=-1/2$ for $(M,\mathcal{K})=(\imath{\rm o}(2n),{\rm O}(2n))$ and $\nu=+1/2$ for $(M,\mathcal{K})=(\imath{\rm o}(2n+1),{\rm O}(2n+1)),(\imath{\rm usp}(2n),{\rm USp}(2n))$. The embedding $\iota_{M}$ is for the single matrix spaces
 			\begin{equation}\label{embedding}
 			\begin{split}
 			\iota_{{\rm O}(2n)}(a)=&\sqrt{a}\otimes\tau_2,\ \iota_{{\rm O}(2n+1)}(a)=\diag(\sqrt{a}\otimes\tau_2,0),\ \iota_{K_4}(a)=\sqrt{a}\otimes\tau_3,\ \iota_{M_\nu}(a)=\left[\begin{array}{cc} 0 & \sqrt{a}\Pi_{n,n+\nu} \\ \Pi_{n,n+\nu}^*\sqrt{a} & 0 \end{array}\right],
 			\end{split}
 			\end{equation}
 			where $\Pi_{ab}$ is the projection from $b$ rows onto the first $a$ rows and $\tau_2$ the second Pauli matrix.
Note we do not distinguish between the two cases $(\imath{\rm o}(2n+1),{\rm O}(2n+1))$ and $(\imath{\rm usp}(2n),{\rm USp}(2n))$ for  $\nu=+1/2$ since the spectral statistics are exactly the same~\cite[Lemma~3.4]{Kieburg:2017} for $\mathcal{K}$-invariant random matrix ensembles.
We employed the constants
\begin{equation}\label{const}
 C_n= \frac{1}{n!}\prod_{j=0}^{n-1}\frac{\pi^j}{j!}\qquad {\rm and}\qquad  C_{n,\nu}^*= \frac{1}{n!}\prod_{j=0}^{n-1}\frac{\pi^{2j+\nu+1}}{\Gamma[j+\nu+1]j!}
\end{equation}
with $\Gamma$ being the Gamma function and we used the following convention for the Vandermonde determinant,
\begin{equation}\label{Vandermonde}
\Delta_n(a)=\prod_{1\leq b<c\leq n}(a_c-a_b)=\det[a_l^{k-1}]_{l,k=1,\ldots,n}.
\end{equation}
Moreover, the subsets of probability densities of these sets will be denoted by the subscript ``Prob", e.g. $L_{\rm Prob}^{1,K_2}(H_2)$.

Our major interest lies in the convolutions on $M=H_2$ and on $M=H_1,H_4,M_\nu$ which are given by
 			\begin{equation}\label{conv-def}
 			f_{H_2}\ast h_{H_2}(y)=\int_{H_2}f_{H_2}(y')h_{H_2}(y-y')dy'\quad {\rm and}\quad f_{M}\ast_\nu h_{M}(y)=\int_{M}f_{M}(y')h_{M}(y-y')dy'
 			\end{equation}
for any two functions $f_{M},h_{M}\in L^{1,\mathcal{K}}(M)$ with $\mathcal{K}=K_\beta,\hat{K}_\nu$, respectively. For this purpose we concentrate on polynomial ensembles and their subsets called P\'olya ensembles. We want to briefly recall their definitions. To do this we need the following subsets of $L^1$-functions
\begin{equation}\label{L1-sets}
 \begin{split}
 L^{1}_{[1,n]}(\mathbb{R})=&\biggl\{f\in L^1(\mathbb{R})\biggl| \text{for all $\kappa\in[1,n]$}:\ \int_{-\infty}^\infty\left|x^{\kappa-1}f(x)\right|dx<\infty\biggl\},\\
 L^{1}_{\mathcal{F}}(\mathbb{R})=&\biggl\{f\in L^1(\mathbb{R})\biggl| f\text{ is non-negative and $(n-1)$-times differentiable and}\\
 &\text{for all $\kappa\in[1,n]$ and }j=0,\ldots,n-1:\ \int_{-\infty}^\infty\left|x^{\kappa-1}\frac{\partial^jf}{\partial x^j}(x)\right|dx<\infty\biggl\},\\
 L^{1}_{\nu}(\mathbb{R}_+)=&\biggl\{f\in L^1(\mathbb{R}_+)\biggl| f\text{ is  non-negative and $2(n-1)$-times differentiable,}\\
 &\hspace*{-1.5cm}\text{for all $\kappa\in[1,n]$ and }j=0,\ldots,n-1:\ \int_{0}^\infty\left|x^{\kappa-1}\left(x^{\nu}\frac{\partial}{\partial x}\frac{1}{x^{\nu-1}}\frac{\partial}{\partial x}\right)^jf(x)\right|dx<\infty,\\
 &\hspace*{-1.5cm}\text{and }\lim_{x\to0} x^{\nu+1}\frac{\partial}{\partial x}\frac{1}{x^\nu}\left(\frac{\partial}{\partial x}x^{\nu+1}\frac{\partial}{\partial x}\frac{1}{x^\nu}\right)^lf(x)=0\ \text{for all }l=0,\ldots,n-2\biggl\}.
 \end{split}
\end{equation}

\begin{definition}[Polynomial and P\'olya ensembles]\label{def:PolEns} \

\begin{enumerate}
\item
A probability density  $p_D\in L_{\rm Prob}^{1,\mathbb{S}}(\mathcal{D})$ is called the \emph{polynomial ensemble} on $\mathcal{D}=A,D$
associated with the one-point weights $w_1,\ldots,w_n\in L^{1}_{[1,n]}(R)$ with $R=\mathbb{R},\mathbb{R}_+$, respectively, 
if it has the form~\cite{Kuijlaars:2014}
\begin{equation}\label{jpdf-pol-ens}
	p_D(a)=\frac{C_n[w]}{n!}\Delta_n(a)\det[w_b(a_c)]_{b,c=1,\ldots,n} \ge 0,\quad a \in \mathcal{D},
\end{equation}
with $C_n[w] > 0$ the normalization constant. 
\item
A probability measure $p_{M}\in L_{\rm Prob}^{1,\mathcal{K}}(M)$ with $M=H_\beta,M_\nu$ and $\mathcal{K}=K_\beta,\hat{K}_\nu$ is called a polynomial ensemble on $M$
if the corresponding eigenvalue (squared singular value) distribution is a polynomial ensemble on $\mathcal{D}$.
\item A polynomial ensemble on $M$ is called P\'olya ensemble on~$M$ iff (\cite[Sec.~3.4]{Kuijlaars:2016a} and~\cite[Definition~3.5]{Kieburg:2017})
\begin{equation}\label{pol-der-1}
w_j(x)=\left(-\frac{\partial}{\partial x}\right)^{j-1}\omega(x),\ \text{ for all }x\in\mathbb{R} \text{ and }j=1,\ldots,n\ {\rm 
with}\ \omega\in  L^{1}_{\mathcal{F}}(\mathbb{R}),
\end{equation}
for $M=H_2$ or~\cite[Definition~3.5]{Kieburg:2017}
\begin{equation}\label{pol-der-2}
w_j(x)=\left(x^\nu\frac{\partial}{\partial x}x^{1-\nu}\frac{\partial}{\partial x}\right)^{j-1}\omega(x),\ \text{ for all }x\in\mathbb{R}_+ \text{ and }j=1,\ldots,n\ {\rm 
with}\ \omega\in  L^{1}_{\nu}(\mathbb{R}_+)
\end{equation}
for $M=H_1,H_4,M_\nu$.
\end{enumerate}
\end{definition}

It was shown in~\cite{Kuijlaars:2016a,Kieburg:2017} that the matrix convolution of a P\'olya ensemble with a polynomial ensemble on the same set of matrices yields again a polynomial ensemble. Additionally the matrix convolution of two P\'olya ensembles of the same kind is closed. The reason for this is that the convolution on the matrix level can be traced back to the additive convolution on $\mathbb{R}$ for $M=H_2$ and to the additive convolution of radially symmetric functions on $\mathbb{R}^{2\nu+2}$ for $M=H_1,H_4,M_\nu$. We denote these convolutions also with ``$\ast$" and ``$\ast_\nu$" since they are related to the convolution~\eqref{conv-def}. In particular the convolution on $\mathbb{R}$ is
\begin{equation}\label{uni-conv-1}
\omega\ast\sigma(x)=\int_{-\infty}^\infty \omega(y)\sigma(x-y)dy
\end{equation} 
while the convolution on $\mathbb{R}^{2\nu+2}$ reduced to the radial part is
\begin{equation}\label{uni-conv-2}
\begin{split}
\omega\ast_\nu\sigma(x)=&x^\nu\int_{0}^\infty \omega(y)\left(\int_{{\rm O}(2\nu+2)/{\rm O}(2\nu+1)}\frac{\sigma(||\sqrt{x}e_1-\sqrt{y}e_\varphi||^2)}{||\sqrt{x}e_1-\sqrt{y}e_\varphi||^{2\nu}}d^*e_\varphi\right)dy\\
 			=&\left\{\begin{array}{cl} \displaystyle\frac{\Gamma[\nu+1]}{\sqrt{\pi}\Gamma[\nu+1/2]}x^\nu\int_{0}^\infty \omega(y)\left(\int_{-1}^1\frac{\sigma(x+y-2\sqrt{xy} t)}{(x+y-2\sqrt{xy} t)^{\nu}}(1-t^2)^{\nu-1/2}dt\right)dy, & {\rm for}\ \nu>-1/2,\\ \displaystyle\frac{1}{2\sqrt{x}}\int_{0}^\infty \omega(y)\bigl[|\sqrt{y}-\sqrt{x}|\sigma((\sqrt{y}-\sqrt{x})^2)+|\sqrt{y}+\sqrt{x}|\sigma((\sqrt{y}+\sqrt{x})^2)\bigl]dy , & {\rm for}\ \nu=-1/2 \end{array}\right.
\end{split}
\end{equation}
for two suitably integrable functions $\omega$ and $\sigma$.
The vector $e_1\in\mathbb{R}^{2\nu+2}$ is some fixed unit vector  while the vector $e_\varphi\in\mathbb{R}^{2\nu+2}$ parametrizes the $2\nu+1$ dimensional unit sphere ${\rm O}(2\nu+2)/{\rm O}(2\nu+1)$ (${\rm O}(0)=1$) and $||.||$ is the Euclidean norm. We underline that we consider the convolutions on the level of densities which explains the factors $x^\nu$ and $||\sqrt{x}e_1-\sqrt{y}ke_1||^{-2\nu}$ in Eq.~\eqref{uni-conv-2}. On the level of functions we have to convolute $\omega(||v||^2)/||v||^{2\nu}$ and $\sigma(||v||^2)/||v||^{2\nu}$ with $v\in\mathbb{R}^{2\nu+2}$. The convolutions~\eqref{uni-conv-1} and \eqref{uni-conv-2} are related to the following two univariate transforms. The first transform is the Fourier transform of $f\in L^1(\mathbb{R})$,
 			\begin{equation}\label{F-def}
 			\begin{split}
 			\mathcal{F}f(s)=&\int_{-\infty}^\infty f(x)\exp[\imath xs]dx,
 			\end{split}
 			\end{equation}
which is found for the additive convolution on $H_2$. The second transform is the ``modified" Hankel transform of $f\in L^1(\mathbb{R}_+)$ which is
 			\begin{equation}\label{H-def}
 			\begin{split}
 			\mathcal{H}_\nu f(s)=&\Gamma[\nu+1]\int_{0}^\infty f(x) \frac{J_\nu(2\sqrt{xs})}{(sx)^{\nu/2}} dx
 			\end{split}
 			\end{equation}
with $\nu\in\mathbb{R}$. This transform naturally appears for the case $M=H_1,H_4,M_\nu$, see~\cite{Kieburg:2017}. Indeed the relation between the transforms and the convolutions are
\begin{equation}
\mathcal{F}[\omega\ast\sigma]=\mathcal{F}\omega\mathcal{F}\sigma \quad {\rm and}\quad \mathcal{H}_\nu[\omega\ast_\nu\sigma]=\mathcal{H}_\nu\omega\mathcal{H}_\nu\sigma,
\end{equation}
respectively. With help of these transforms one can easily derive the explicit normalization constants for the two kinds of P\'olya ensembles.

\begin{lemma}[JPDF of P\'olya Ensembles]\label{lem:jpdf.Polya}\
\begin{enumerate}
\item	Let $X\in H_2$ be a random matrix drawn from the P\'olya ensemble on $H_2$ associated to the weight $\omega\in L^{1}_{\mathcal{F}}(\mathbb{R})$. The joint probability density of the unordered eigenvalues $x\in D$ of $X$ is given by
			\begin{equation}\label{jpdf.Polya.a}
			p_D(x)=\frac{1}{n!}\left(\prod_{j=0}^{n-1}\frac{1}{\mathcal{F}\omega(0) j!}\right)\Delta_n(x)\det\left[(-\partial_a)^{b-1}\omega(x_a)\right]_{a,b=1,\ldots,n}.
			\end{equation}
\item	Let $X\in M=H_1,H_4,M_\nu$ be a random matrix drawn from the P\'olya ensemble on $M$ associated to the weight $\omega\in L^{1}_{\nu}(\mathbb{R}_+)$. The joint probability density of the unordered squared singular values $x\in A$ of $X$ is given by
			\begin{equation}\label{jpdf.Polya.b}
			p_A(x)=\frac{1}{n!}\left(\prod_{j=0}^{n-1}\frac{\Gamma[\nu+1]}{\mathcal{H}_\nu\omega(0) j!\Gamma[\nu+j+1]}\right)\Delta_n(x)\det\left[(\partial_a x_a^{\nu+1}\partial_a x_a^{-\nu})^{b-1}\omega(x_a)\right]_{a,b=1,\ldots,n}.
			\end{equation}
\end{enumerate}
\end{lemma}

{\bf Proof:}\\
First we want to underline that $x_a^\nu\partial_ax_a^{1-\nu}\partial_a=\partial_a x_a^{\nu+1}\partial_a x_a^{-\nu}$. To prove both statements in a unifying way, we integrate over all eigenvalues/squared singular values $x$ and apply Andr\'eief's integration theorem~\cite{Andreief:1883}. This reduces the problem to one-dimensional integrals where we can integrate by parts to apply the derivatives onto the monomials resulting from the Vandermonde determinant. Since both operators  $\partial_a$ and $x_a^{-\nu}\partial_a x_a^{\nu+1}\partial_a $ map monomials of order $j$ to those of order $j-1$ and annihilate constants we have to take the determinant of an upper triangular matrix. The diagonal elements are exactly the factors in the products~\eqref{jpdf.Polya.a} and~\eqref{jpdf.Polya.b}, respectively.\hfill$\square$
\vskip0.2cm

An important property of P\'olya ensembles are particular simple group integrals, see \cite[Theorem 4.6]{Kieburg:2017}. Those integrals can be employed to calculate the joint probability density of a P\'olya ensemble shifted by a fixed matrix.

\begin{theorem}[JPDF of P\'olya Ensembles Convoluted with Fixed Matrices]\label{thm:jpdf.Polya.fixed}\

\begin{enumerate}
\item 	 Let $X_1\in H_2$ be a random matrix drawn from the P\'olya ensemble on $H_2$ associated to the weight $\omega\in L^{1}_{\mathcal{F}}(\mathbb{R})$ and $X_2\in H_2$ be a fixed matrix with non-degenerate eigenvalues $x\in D$. The joint probability density of the eigenvalues $y\in D$ of the random matrix $Y=X_1+X_2\in H_2$ is given by
			\begin{equation}\label{jpdf.fixed.a}
 			\begin{split}
 			p_D(y|x)=&\frac{1}{n!(\mathcal{F}\omega(0))^n}\frac{\Delta_n(y)}{\Delta_n(x)}\det[\omega(y_a-x_b)]_{a,b=1,\ldots,n},
 			\end{split}
 			\end{equation}
 			i.e. it is a polynomial ensemble associated to the weights $\{\omega(.-x_j)\}_{j=1,\ldots,n}$. 
\item	 Let $X_1\in M=H_1,H_4,M_\nu$ be a random matrix drawn from the P\'olya ensemble on $M$ associated to the weight $\omega\in L^{1}_{\nu}(\mathbb{R}_+)$ and $X_2\in M$ be a fixed matrix with non-degenerate eigenvalues $x\in A$. The joint probability density of the eigenvalues $y\in A$ of the random matrix $Y=X_1+X_2\in M$ is given by
			\begin{equation}\label{jpdf.fixed.b}
 			\begin{split}
 			p_A(y|x)=&\frac{1}{n!}\left(\prod_{j=0}^{n-1}\frac{\Gamma[\nu+1]}{\sqrt{\pi}\Gamma[\nu+1/2]\mathcal{H}_\nu\omega(0)}\right)\frac{\Delta_n(y)}{\Delta_n(x)}\det\left[y_a^\nu\int_{-1}^{1} \frac{\omega(y_a+x_b-2\sqrt{y_ax_b}t)}{(y_a+x_b-2\sqrt{y_ax_b}t)^{\nu}}(1-t^2)^{\nu-1/2}dt\right]_{a,b=1,\ldots,n}
 			\end{split}
 			\end{equation}
 			for $\nu>-1/2$ and
			\begin{equation}\label{jpdf.fixed.c}
 			\begin{split}
 			p_A(y|x)=&\frac{1}{n!(\mathcal{H}_{-1/2}\omega(0))^n}\frac{\Delta_n(y)}{\Delta_n(x)}\det\left[\frac{|\sqrt{y_a}-\sqrt{x_b}|\omega((\sqrt{y_a}-\sqrt{x_b})^2)+|\sqrt{y_a}+\sqrt{x_b}|\omega((\sqrt{y_a}+\sqrt{x_b})^2)}{2\sqrt{y_a}}\right]_{a,b=1,\ldots,n}
 			\end{split}
 			\end{equation}
 			for $\nu=-1/2$.
 			Thus it is again a polynomial ensemble.
\end{enumerate}
\end{theorem}

{\bf Proof:}\\
Again we prove both statements at the same time. Let $P_M(X_1)$ be the distribution of the P\'olya ensemble on $M$ corresponding to $\omega$. Hence the joint probability distribution of the eigenvalues/squared singular values $y$ of $Y=X_1+X_2$ is,  up to a constant, given  by
\begin{equation}
p_{\mathcal{D}}(y|x)\propto \det y^\nu \Delta_n^2(y) \int_{\mathcal{K}} P_M(X_2+k\iota_M(y)k^*)d^*k
\end{equation}
with $\iota_M$ as in Eq.~\eqref{embedding} and $\mathcal{D}=D,A$ and $\mathcal{K}=K_\beta,\hat{K}_\nu$, respectively. The matrix $X_2$ can be also decomposed as $X_2=\tilde{k}\iota_M(x)\tilde{k}^*$ with $x\in\mathcal{D}$ and $\tilde{k}\in\mathcal{K}$. Due to the $\mathcal{K}$-invariance of $P_M$ we have $P_M(\tilde{k}\iota_M(x)\tilde{k}^*+k\iota_M(y)k^*)=P_M(\iota_M(x)+\tilde{k}^*k\iota_M(y)k^*\tilde{k})$ and we can absorb  the unitary matrix $\tilde{k}$ in the group integral. The group integral was calculated for any P\'olya ensemble in \cite[Theorem 4.6]{Kieburg:2017} and is
\begin{equation}
\int_{\mathcal{K}} P_M(X_2+k\iota_M(y)k^*)d^*k\propto\frac{1}{\Delta_n(x)\Delta_n(y)}\det\left[\int_{\mathcal{K}|_{n=1}} \{P_M(\iota_M(x_b)+k\iota_M(y_a)k^*)\}|_{n=1}d^*k\right]_{a,b=1,\ldots,n},
\end{equation}
where $\{.\}|_{n=1}$ means that it is the distribution and integral for the matrix spaces with the dimension parameter $n=1$. More explicitly we have
\begin{equation}
\{P_M(\iota_M(x_b)+k\iota_M(y_a)k^*)\}|_{n=1}\propto\frac{\omega(x_b+y_a)}{\mathcal{F}\omega(0)}
\end{equation}
for $M=H_2 $ and
\begin{equation}
\{P_M(\iota_M(x_b)+k\iota_M(y_a)k^*)\}|_{n=1}\propto\frac{\omega(\tr(\iota_M(x_b)+k\iota_M(y_a)k^*)^2/2)}{\mathcal{H}_\nu\omega(0)(\tr(\iota_M(x_b)+k\iota_M(y_a)k^*)^2/2)^{\nu}}
\end{equation}
for the other cases, because of the immersion $\mathcal{I}_M$, see Eq.~\eqref{I-H2}. While for $M=H_2$ the group integral drops out we have a remaining integral for the other cases over a $2\nu+1$ dimensional unit sphere. This can be seen by noticing
\begin{equation}
\begin{split}
\tr(\iota_{\imath{\rm o}(2)}(x_b)k\iota_{\imath{\rm o}(2)}(y_a)k^*)=&2\det(k)\sqrt{x_by_a},\\
\tr(\iota_{\imath{\rm o}(3)}(x_b)k\iota_{\imath{\rm o}(3)}(y_a)k^*)=&2\sqrt{x_by_a} (k_{11}k_{22}-k_{12}k_{21}),\\
\tr(\iota_{\imath{\rm usp}(2)}(x_b)k\iota_{\imath{\rm usp}(2)}(y_a)k^*)=&2\sqrt{x_by_a}\tilde{k}_{33},\\
\tr(\iota_{M_\nu}(x_b)k\iota_{M_\nu}(y_a)k^*)=&2\sqrt{x_by_a} {\rm Re}(k_{11}k_{\nu+2,\nu+2}^*),
\end{split}
\end{equation}
where we used the relation
\begin{equation}
\begin{split}
kk^T=&\eins_3\ \rightarrow\ |k_{11}k_{22}-k_{12}k_{21}|=|k_{33}|,\qquad {\rm for}\ M=\imath{\rm o}(3),
\end{split}
\end{equation}
and that the adjoint representation $\{\tr k\tau_a k^* \tau_b/2\}_{a,b=1,2,3}=\{\tilde k\}_{a,b=1,2,3}\in{\rm SO}(3)$ of $k\in{\rm USp}(2)$ is the three-dimensional special orthogonal group also distributed by the Haar measure.
Note that for $M=\imath{\rm o}(3),\imath{\rm usp}(2)$ the third row of $k$ parametrizes a two-dimensional sphere and that for $M=M_\nu$ the first row of $k\in \hat{K}_\nu$ is given by a $2\nu+1$ dimensional unit sphere. The integral over the unit sphere only depends on a single component of the sphere which is given by the parametrization in Eq.~\eqref{jpdf.fixed.b}.  Moreover we equip the zero-dimensional unit sphere, which is $\mathbb{Z}_2$, with the normalized Dirac measure at its two elements.

The normalization constant for $M=H_2$ can be found in the limit $x\to0$ which has to yield the result~\eqref{jpdf.Polya.a}. For the other case of $M$ we can readily fix the normalization constant by the particular choice of the Laguerre ensemble, i.e. $\omega(z)=z^\nu e^{-z}$, in the limit $x\to0$ since the constant is independent of $\omega$ and $x$. The integral over $t$ in Eq.~\eqref{jpdf.fixed.b} yields the renormalized modified Bessel function of the first kind $I_\nu(2\sqrt{z})/z^{\nu/2}$ and  for the limit $x\to0$ we have to apply l'H\^{o}spital's rule yielding the prefactor of Eq.~\eqref{jpdf.fixed.b}. This closes the proof.
\hfill$\square$
\vskip0.2cm

The exact statement for the matrix convolutions with P\'olya ensembles is closely related to Theorem~\ref{thm:jpdf.Polya.fixed}.

\begin{theorem}[JPDF of P\'olya Ensembles Convoluted with Polynomial Ensembles]\label{thm:jpdf.Polya.poly}\

\begin{enumerate}
\item 	 Let $X_1\in H_2$ be a random matrix drawn from the P\'olya ensemble on $H_2$ associated to the weight $\omega\in L^{1}_{\mathcal{F}}(\mathbb{R})$ and $X_2\in H_2$ be a random matrix drawn from a polynomial ensemble on $H_2$ associated with the weights $w_1,\hdots,w_n$. The joint probability density of the eigenvalues $y\in D$ of the random matrix $Y=X_1+X_2\in H_2$ is given by
			\begin{equation}\label{jpdf.poly.a}
 			\begin{split}
 			p_D^{(w)}(y)=&\frac{C_n[w]}{n!(\mathcal{F}\omega(0))^n}\Delta_n(y)\det[\omega\ast w_b(y_a)]_{a,b=1,\ldots,n}
 			\end{split}
 			\end{equation}
 			and, thus, is again a polynomial ensemble associated to the weights $\{\omega\ast w_j\}_{j=1,\ldots,n}$. In the case that $X_2$ is also drawn from a P\'olya ensemble on $H_2$ associated with the weight $\sigma$ the random matrix $Y$ is a P\'olya ensemble on $H_2$ associated with the weight $\omega\ast\sigma$.
			(This statement was proven apart from the normalization in~\cite{Kuijlaars:2016a}.)
\item	 Let $X_1\in M=H_1,H_4,M_\nu$ be a random matrix drawn from the P\'olya ensemble on $M$ associated to the weight $\omega\in L^{1}_{\nu}(\mathbb{R}_+)$ and $X_2\in M$be a random matrix drawn from a polynomial ensemble on $M$ associated with the weights $w_1,\hdots,w_n$. The joint probability density of the eigenvalues $y\in A$ of the random matrix $Y=X_1+X_2\in M$ is
			\begin{equation}\label{jpdf.poly.b}
 			\begin{split}
 			p_A^{(w)}(y)=&\frac{C_n[w]}{n!(\mathcal{H}_\nu\omega(0))^n}\Delta_n(y)\det\left[\omega\ast_\nu w_b(y_a)\right]_{a,b=1,\ldots,n},
 			\end{split}
 			\end{equation}
 			which is a polynomial ensemble associated with the weights $\{\omega\ast_\nu w_j\}_{j=1,\ldots,n}$.  If $X_2$ is  drawn from a P\'olya ensemble on $M$ associated with the weight $\sigma$, too, the random matrix $Y$ is a P\'olya ensemble on $M$ with the weight $\omega\ast_\nu\sigma$.
			(This statement was proven in~\cite{Kieburg:2017}.)
\end{enumerate}
\end{theorem}

{\bf Proof:}\\
The proofs of the two statements were essentially done  in~\cite[Corollary~3.3. and 3.4]{Kuijlaars:2016a} and in~\cite[Theorem 3.10 and Corollary 3.11]{Kieburg:2017}. The normalization can be read off from Theorem~\ref{thm:jpdf.Polya.fixed} and the Definition~\eqref{jpdf-pol-ens}. One has only to integrate over $x\in \mathcal{D}=D,A$ weighted by the distribution~\eqref{jpdf-pol-ens} which are the eigenvalues/squared singular values of $X_2$. The Vandermonde determinant $\Delta_n(x)$ cancels and the integral can be done by Andr\'eief's identity~\cite{Andreief:1883}.
\hfill$\square$
\vskip0.2cm

The aim of the present work is to construct the bi-orthogonal functions and kernels for the three cases of joint probability densities presented in Lemma~\ref{lem:jpdf.Polya}, Theorem~\ref{thm:jpdf.Polya.fixed} and Theorem~\ref{thm:jpdf.Polya.poly}. For this reason we call the pair of functions $\{p_j,q_j\}_{j=0,\ldots,n-1}$ a pair of bi-orthonormal functions on $R=\mathbb{R},\mathbb{R}_+$ when it satisfies
\begin{equation}
 \int_R p_l(x)q_m(x) dx=\delta_{lm},\ {\rm for\ all\ }l,m=0,\ldots,n-1,
\end{equation}
with $\delta_{lm}$ the Kronecker symbol. Then a polynomial ensemble can be described by the pair of bi-orthonormal functions $\{p_j,q_j\}_{j=0,\ldots,n-1}$ if the linear spans of $\{p_j\}_{j=0,\ldots,n-1}$ and $\{q_j\}_{j=0,\ldots,n-1}$ are equal to the linear spans of $\{x^j\}_{j=0,\ldots,n-1}$ and $\{w_{j+1}\}_{j=0,\ldots,n-1}$, respectively. In particular the $k$-point correlation function of the polynomial ensemble~\eqref{jpdf-pol-ens} has the form~\cite{Borodin:1999}
\begin{equation}\label{k-point}
R_k(x_1,\ldots,x_k)=\det[K_n(x_b,x_c)]_{b,c=1,\ldots,n}\ {\rm with}\ K_n(x_b,x_c)=\sum_{j=0}^{n-1}p_j(x_b)q_{j}(x_c).
\end{equation}
Thus the whole statistics are determined when a pair of bi-orthonormal functions of the polynomial ensemble is known.

\section{P\'olya Ensembles on $H_2$}\label{sec:H2}

We first consider the structurally simpler case of P\'olya ensembles on $H_2$. In subsection~\ref{sec:evstat.Polya.a} we derive the bi-orthonormal functions of a general P\'olya ensemble on $H_2$ without any shift. In the same section we point out a relation to Toeplitz determinants of which the author is not aware that it already exists in the literature in this generality as shown. The case of a shift by a fixed matrix in $H_2$ and by a polynomial ensemble on $H_2$ are considered in subsections~\ref{sec:evstat.fixed.a} and~\ref{sec:evstat.poly.a}, respectively.

\subsection{Statistics of P\'olya Ensembles on $H_2$}\label{sec:evstat.Polya.a}

We will first state and prove a theorem which applies for any P\'olya ensemble on $H_2$. Later we are going to rephrase it into a much simpler and more recognizable form when the weight $\omega$ satisfies additional analyticity properties. 

\begin{theorem}[Eigenvalue Statistics of P\'olya Ensembles]\label{thm:stat.Polya.a}\

 The joint probability density~\eqref{jpdf.Polya.a} can be described by the bi-orthonormal functions
 \begin{equation}\label{biorth.Polya.a}
 \{p_j,q_j\}_{j=0,\ldots,n-1}=\left\{\left.\frac{(y'-\imath\partial_t)^j}{j!}\frac{1}{\mathcal{F}\omega(t)}\right|_{t=0}\,,\,(-\partial_y)^j\omega(y)\right\}_{j=0,\ldots,n-1}.
 \end{equation}
 The kernel is given by
 \begin{equation}\label{kernel.Polya.a}
 \begin{split}
  K_n(y',y)=&\int_{0}^\infty ds e^{-s}\left.(s-(y'-\imath\partial_t)\partial_y)^{n-1}\frac{\omega(y)}{\mathcal{F}\omega(t)}\right|_{t=0}.
  \end{split}
 \end{equation}
\end{theorem}

{\bf Proof:}\\
First we prove the bi-orthonormality of the pair of functions which is the integral
\begin{equation}
I_{lm}=(-\imath)^{l}\int_{-\infty}^\infty dy\left.\frac{(y-\imath\partial_t)^l}{l!}\frac{1}{\mathcal{F}\omega(t)}\right|_{t=0}(-\partial_y)^m\omega(y).
\end{equation}
In the first step we integrate  by parts in $y$. The boundary terms vanish due to the integrability and differentiability conditions of $\omega\in L_{\mathcal{F}}^1(\mathbb{R})$, cf. Eq.~\eqref{L1-sets}. When doing so we notice that for $l<m$ the integral $I_{lm}$ vanishes because the polynomial is of order $l$ while we differentiate $m$-times. Hence we can reduce the discussion to the case $l\geq m$ and the integral is
\begin{equation}
I_{lm}=\int_{-\infty}^\infty dy\left.\frac{(y-\imath\partial_t)^{l-m}}{(l-m)!}\frac{1}{\mathcal{F}\omega(t)}\right|_{t=0}\omega(y).
\end{equation}
In the next step we employ the identity
\begin{equation}
\int_{-\infty}^\infty dy y^a \omega(y)=\left.(-\imath\partial_s)^a\mathcal{F}\omega(s)\right|_{s=0},\ {\rm for\ all}\ a=0,\ldots,n-1,
\end{equation}
which is again true because of the integrability and differentiability conditions of $\omega$. This yields
\begin{equation}
I_{lm}=\left.\frac{\left(-\imath\partial_s-\imath\partial_t\right)^{l-m}}{(l-m)!}\frac{\mathcal{F}\omega(s)}{\mathcal{F}\omega(t)}\right|_{s=t=0}.
\end{equation}
Changing to relative, $s-t$, and center of mass, $s+t$, coordinates it becomes immediate that for $l>m$ the derivative vanishes while for $l=m$ we have unity. This proves the bi-orthonormality.

For deriving the kernel~\eqref{kernel.Polya.a} we start from the standard form~\eqref{k-point}
\begin{equation}
\begin{split}
K_n(y',y)=&\sum_{j=0}^{n-1} \left.\frac{(y'-\imath\partial_t)^j}{j!}\frac{1}{\mathcal{F}\omega(t)}\right|_{t=0}(-\partial_y)^j\omega(y)=\int_{0}^\infty ds e^{-s}\sum_{j=0}^{n-1} \left.\frac{(y'-\imath\partial_t)^js^{n-1-j}}{j!(n-1-j)!}\frac{1}{\mathcal{F}\omega(t)}\right|_{t=0}(-\partial_y)^j\omega(y).
\end{split}
\end{equation}
The sum is the binomial sum yielding the claim.
\hfill$\square$
\vskip0.2cm

We want to point out one particular result which can be immediately derived from Theorem~\ref{thm:stat.Polya.a} which relates two Toeplitz determinants. For this purpose we want to consider the average
\begin{equation}
Z_L(z)=\int_D \prod_{a=1}^L\prod_{b=1}^n(z_a-y_b) p_D(y)dy
\end{equation}
for $L\leq n-1$ being a positive integer and $z_1,\ldots,z_L\in\mathbb{C}$ being pairwise different complex numbers. We underline that we have to modify the integrability conditions of $\omega\in L_{\mathcal{F}}^1(\mathbb{R})$ in Eq.~\eqref{L1-sets} to guarantee the integrability of this average.  The product in front of the joint probability density can be combined with the Vandermonde determinant in Eq.~\eqref{jpdf.Polya.a} and we can perform a generalized version of Andr\'eief's identity~\cite[Appendix~C.1]{Kieburg:2010} which yields
\begin{equation}
Z_L(z)=\frac{\prod_{j=n}^{n+L-1}j!}{\Delta_L(z)}\det\left[\begin{array}{c} \displaystyle\Theta_{ac}\left.\frac{(-\imath\partial_t)^{c-a}}{(c-a)!}\frac{\mathcal{F}\omega(t)}{\mathcal{F}\omega(0)}\right|_{t=0} \\ \displaystyle\frac{z_b^{c-1}}{(c-1)!} \end{array}\right]_{\substack{ a=1,\ldots,n \\ b=1,\ldots,L \\ c=1,\ldots,n+L}},
\end{equation}
where $\Theta_{ac}=1$ for $a\leq c$ and otherwise vanishes.  Moreover we know the bi-orthogonal polynomials corresponding to this ensemble and so
\begin{equation}
Z_L(z)=\frac{1}{\Delta_L(z)}\det\left[\begin{array}{c} (n+b-1)!\mathcal{F}\omega(0)p_{n+b-1}(z_a) \end{array}\right]_{\substack{ a,b=1,\ldots,L}},
\end{equation}
see \cite{Uvarov:1969}.
The constants in the product in front of the polynomials correctly normalize them to monic normalization. In the next step we take the limit $z\to0$ and find the following identity between the two Toeplitz determinants
 \begin{equation}\label{toeplitz}
 \det\left[\begin{array}{ccccccc} c_0 & c_1 & & \cdots &  & c_{n-2} & c_{n-1} \\  c_{-1} & c_0 & & \cdots & & & c_{n-2}\\ \vdots &  &  & \\ c_{-L} & & & \ddots & & \vdots & \vdots \\ 0 & & & \\ \vdots & \ddots & &  &  & c_0 & c_1 \\ 0 &\cdots & 0 & c_{-L} & \cdots & c_{-1} & c_0 \end{array}\right]=(-1)^{nL}\det\left[\begin{array}{ccccc} d_{L-1} & d_{L} & \cdots & d_{2L-2} & d_{2L-1} \\ d_{L-2} & d_{L-1} & \cdots &  & d_{2L-2} \\ \vdots & \vdots & \ddots & \vdots & \vdots \\  d_{1} &  & \cdots & d_{L-1} & d_{L} \\  d_{0} & d_1 & \cdots & d_{1} & d_{L-1} \end{array}\right]
 \end{equation}
 with
 \begin{equation}\label{coefficients}
 c_j=\left.\frac{(-\imath\partial_t)^{L+j}}{(L+j)!}\frac{\mathcal{F}\omega(t)}{\mathcal{F}\omega(0)}\right|_{t=0}\quad{\rm and}\quad d_{L-1+b-a}=\left.\frac{(-\imath\partial_t)^{n+b-a}}{(n+b-a)!}\frac{\mathcal{F}\omega(0)}{\mathcal{F}\omega(t)}\right|_{t=0}.
 \end{equation}
 This relation can be generalized to arbitrary Toeplitz determinants of the form~\eqref{toeplitz} since the derivatives of the Fourier transform $\mathcal{F}\omega$ can be quite arbitrary.

\begin{corollary}[Relation between Toeplitz Determinants]\label{cor:toeplitz}\

 Let $L\leq n-1$ be a positive integer and $c_{-L},c_{1-L},\ldots,c_{n-1}\in\mathbb{C}$ arbitrary complex numbers apart from $c_{-L}=1$. We define $F(t)=\sum_{j=0}^{n+L-1}c_{j-L} t^j$ and $d_{L-1+j}=1/(n+j)!\partial_t^{n+j}1/F(t)|_{t=0}$ which replaces the definition~\eqref{coefficients}. Then the relation~\eqref{toeplitz} of the Toeplitz determinants still holds.
\end{corollary}

The case $L=1$ is well-known since it relates the elementary polynomials with $n$ elements with the sum of all homogeneous monomials of a fixed degree, see~\cite[Chapter~4.2]{Prasolov:1994}. The author is not aware that the general form of this statement was derived before. For $L=n-2$ the formula~\eqref{toeplitz} gives a recursion formula from an arbitrary $n\times n$ Toeplitz determinant to an $(n-1)\times(n-1)$ Toeplitz determinant.

{\bf Proof:}\\
We employ the fact
\begin{equation}
c_{b-a}=\left.\frac{\partial_t^{b-a+L}}{(b-a+L)!}F(t)\right|_{t=0}=\left.\frac{\partial_t^{L+b-1}}{(L+b-1)!}t^{a-1}F(t)\right|_{t=0}
\end{equation}
for $b\geq a-L$ because $F(0)=1\neq0$ and all derivatives at $t=0$ exist. For $b< a-L$ the right hand side vanishes. Then we
denote the Toeplitz determinant on the left hand side of Eq.~\eqref{toeplitz} by $T$ and rewrite it as follows
\begin{equation}
T=(-1)^{nL}\lim_{z\to0}\left.\frac{\Delta_{n+L}(\partial_t,z)}{\Delta_L(z)}\Delta_n(t)\prod_{j=1}^n\frac{F(t_j)}{(L+j-1)!}\right|_{t_1=\ldots,t_n=0},
\end{equation}
where we introduced $L$ auxiliary variables $z_1,\ldots,z_L$ which are pairwise different. All derivatives act on everything on the right side. Choosing the polynomials
\begin{equation}
 \tilde{p}_j(x)=\frac{1}{j!}\left.(x+\partial_s)^j\frac{1}{F(s)}\right|_{s=0}
\end{equation}
we rewrite one of the Vandermonde determinants as
\begin{equation}
\Delta_{n+L}(\partial_t,z)=\prod_{j=0}^{n+L-1}j!\ \det\left[\begin{array}{c} \tilde{p}_{c-1}(\partial_{t_a}) \\ \tilde{p}_{c-1}(z_b)  \end{array}\right]_{\substack{ a=1,\ldots,n \\ b=1,\ldots,L \\ c=1,\ldots,n+L}}.
\end{equation}
Since $\tilde{p}_{c-1}(\partial_{t_a})t_a^{b-1}F(t_a)|_{t_a=0}=\delta_{bc}$ we have
\begin{equation}
T=(-1)^{nL}\lim_{z\to0}\frac{\prod_{j=0}^{n-1}j!}{\Delta_L(z)}\det\left[\begin{array}{c} \tilde{p}_{n+c-1}(z_b)  \end{array}\right]_{\substack{b,c=1,\ldots,L }}.
\end{equation}
l'H\^{o}spital's rule yields the claim.
\hfill$\square$
\vskip0.2cm

We want to point out a simplification of the results of Theorem~\ref{thm:stat.Polya.a} when the weight $\omega$ satisfies some additional properties.

\begin{corollary}[Simplification of  Theorem~\ref{thm:stat.Polya.a}]\label{cor:simpl.Polya.a}\

We assume the requirements of  Theorem~\ref{thm:stat.Polya.a} and further assume that the Fourier transform $\mathcal{F}\omega$ is holomorphic at the origin and $z^j\mathcal{F}\omega(z)$ is absolutely integrable along the real line for $j=0,\ldots,n-1$. Then the bi-orthonormal functions can be written as
 \begin{equation}\label{biorth.Polya.a.2}
 \{p_j,q_j\}_{j=0,\ldots,n-1}=\left\{(-\imath)^{j}\oint\frac{dz'}{2\pi\imath {z'}^{j+1}}\frac{\exp[\imath y' z']}{\mathcal{F}\omega(z')}\,,\,\imath^j\int_{-\infty}^\infty \frac{dz}{2\pi} z^j\exp[-\imath yz]\mathcal{F}\omega(z)\right\}_{j=0,\ldots,n-1}.
 \end{equation}
and the kernel has the form
 \begin{equation}\label{kernel.Polya.a.3}
 \begin{split}
  K_n(y',y)=&\oint\frac{dz'}{2\pi}\int_{-\infty}^\infty \frac{dz}{2\pi} \frac{1}{z'-z}\left(1-\left(\frac{z}{z'}\right)^n\right)\exp[\imath(y'z'- yz)]\frac{\mathcal{F}\omega(z)}{\mathcal{F}\omega(z')}.
  \end{split}
 \end{equation}
 The contour of $z'$ encircles the origin $z'=0$ while the contour for $z$ is along the real line.
\end{corollary}

{\bf Proof:}\\
Since the function $z^j\mathcal{F}\omega(z)$ is absolutely integrable for $j=0,\ldots,n-1$ we have
\begin{equation}
(-\partial_y)^{j}\int_{-\infty}^\infty \frac{dz}{2\pi} \exp[-\imath yz]\mathcal{F}\omega(z)=\imath^j\int_{-\infty}^\infty \frac{dz}{2\pi} z^j\exp[-\imath yz]\mathcal{F}\omega(z).
\end{equation}
The holomorphy of  $\mathcal{F}\omega(z)$ and $\mathcal{F}\omega(0)\neq0$ implies that $1/\mathcal{F}\omega(z)$ is also holomorphic at the origin. This allows the calculation
\begin{equation}
\begin{split}
\left.\frac{(y'-\imath\partial_t)^j}{j!}\frac{1}{\mathcal{F}\omega(t)}\right|_{t=0}=&\sum_{j=0}^j\frac{{y'}^{j-l}}{l!(j-l)!}\left.(-\imath\partial_t)^l\frac{1}{\mathcal{F}\omega(t)}\right|_{t=0}=(-\imath)^{j}\oint\frac{dz'}{2\pi\imath {z'}^{j+1}}\frac{\exp[\imath y' z']}{\mathcal{F}\omega(z')},
\end{split}
\end{equation}
where the contour only encircles the origin counter clockwise. The kernel can be easily obtained by doing the geometric sum.
\hfill$\square$
\vskip0.2cm

An important remark is in order. The results of Corollary~\eqref{cor:simpl.Polya.a} resemble results of the supersymmetry method in random matrix theory, e.g. see~\cite{Guhr:2006,Kieburg:2009} without Efetov-Wegner boundary terms (it is the case when we also encircle the point $z=z'$ and without the $1$ in the bracket in Eq.~\eqref{kernel.Polya.a.3}) and~\cite{Kieburg:2011} with the Efetov-Wegner boundary terms. Indeed when identifying the integration variables $z$ and $z'$ with the eigenvalues of a $(1|1)\times(1|1)$ supermatrix $\sigma$, the ratio $z/z'$ is equal to the superdeterminant of $\sigma$ and the term $1/(z'-z)$ is the result of the Berezinian (Jacobian in superspace) after diagonalizing $\sigma$ like the Vandermonde determinant for ordinary matrices and the supergroup integral also resulting from the diagonalization. These results of Corollary~\eqref{cor:simpl.Polya.a} also resemble expressions derived via other methods, e.g. see~\cite{Kuijlaars:2014,Claeys:2015,Kieburg:2015,Kieburg:2016b,Forrester:2017}.

\begin{examples}\label{example.a}\

In~\cite{Kieburg:2017} it was shown that the weight $\omega$ has to be a P\'olya frequency function of order $n$ such that the joint probability density~\eqref{jpdf.Polya.a} belongs to a random matrix ensemble of $n\times n$ Hermitian matrices. P\'olya frequency functions of infinite order have a particular simple and explicit expression in terms of its Laplace transform~\cite{Schoenberg:1951,Karlin:1968} or as we write it in terms of the Fourier transform which is either of the form
\begin{equation}\label{Laplace-Polya.a}
\begin{split}
\mathcal{F}\omega(s)=&\exp[-\gamma s^2]\prod_{j=1}^\infty\frac{\exp[-\imath\delta_j s]}{1-\imath\delta_j s},\\
 \gamma\geq0,\ \delta_j\in\mathbb{R},\ 0<\gamma+&\sum_{j=1}^\infty \delta_j^2<\infty,\ {\rm and}\ -\max_{\delta_j<0}\left\{\frac{1}{\delta_j}\right\}>{\rm Im}\, s>-\min_{\delta_j>0}\left\{\frac{1}{\delta_j}\right\},
 \end{split}
\end{equation}
corresponding to a support of $\omega$ on $\mathbb{R}$ or of the form
\begin{equation}\label{Laplace-Polya.b}
\begin{split}
\mathcal{F}\omega(s)=&\prod_{j=1}^\infty\frac{1}{1-\imath \delta_j s},\quad \delta_j\geq0,\, 0<\sum_{j=1}^\infty \delta_j<\infty,\ {\rm and}\ {\rm Im}\, s>-\min\left\{\frac{1}{\delta_j}\right\},
 \end{split}
\end{equation}
corresponding to a support on $\mathbb{R}_0^+$. We omitted the exponential terms in the original work~\cite{Schoenberg:1951} because they only shift the function along the real axis. We have also not included the term $1/s^k$ in the work~\cite{Karlin:1968} since we need to consider integrable weights $\omega$ which excludes a pole of its Fourier transform at $s=0$.

A weight $\omega$ of the form~\eqref{Laplace-Polya.a} is $(n-1)$-times differentiable if either $\gamma>0$ or at least $n+1$ of the $\delta_j$ are non-zero. This follows from the fact that $s^j\mathcal{F}\omega(s)$ is absolutely integrable for $j=0,\ldots, n-1$.  Moreover, $\mathcal{F}\omega$ is holomorphic in a disk about the origin with a radius smaller than $\min_{j=1,2,\ldots}\{1/|\delta_j|\}$ since the derivative with respect to $s$ exists and the Looman-Menchoff theorem~\cite{Looman:1923,Menchoff:1936} can be used, i.e.
\begin{equation}
\partial_s\mathcal{F}\omega(s)=\mathcal{F}\omega(s)\times\left\{\begin{array}{cl} \displaystyle-2\gamma s-\sum_{j=1}^\infty\frac{\delta_j^2}{1-\imath\delta_j s}, & {\rm for\ Eq.~\eqref{Laplace-Polya.a},} \\ \displaystyle-\imath\sum_{j=1}^\infty\frac{\delta_j}{1-\imath\delta_j s}, & {\rm for\ Eq.~\eqref{Laplace-Polya.b}} \end{array}\right.
\end{equation}
is finite inside this disc due to the conditions on $\delta_j$. Additionally we have $\mathcal{F}\omega(0)=1$ and, hence, $\omega$ is normalized due to our choice of the normalization constant. The holomorphy also implies that the integrability conditions of $\omega$ are also automatically satisfied. Collecting everything we can say that the bi-orthonormal pair is given either  by
 \begin{equation}
 \begin{split}
 \{p_j,q_j\}_{j=0,\ldots,n-1}=&\biggl\{\oint\frac{dz'}{2\pi(\imath z')^{j+1}}e^{\gamma {z'}^2+\imath y' z'}\prod_{l=1}^\infty e^{\imath\delta_l z'}(1-\imath\delta_l z')\,,\,\int_{-\infty}^\infty \frac{dz}{2\pi} (\imath z)^je^{-\gamma z^2-\imath yz}\prod_{l=1}^\infty\frac{\exp[-\imath\delta_l z]}{1-\imath\delta_j z}\biggl\}_{j=0,\ldots,n-1}
 \end{split}
 \end{equation}
 or by
 \begin{equation}
 \begin{split}
 \{p_j,q_j\}_{j=0,\ldots,n-1}=&\biggl\{\oint\frac{dz'}{2\pi(\imath z')^{j+1}}e^{\imath y' z'}\prod_{l=1}^\infty (1-\imath\delta_l z')\,,\,\int_{-\infty}^\infty \frac{dz}{2\pi} (\imath z)^je^{-\imath yz}\prod_{l=1}^\infty\frac{1}{1-\imath\delta_l z}\biggl\}_{j=0,\ldots,n-1},
 \end{split}
 \end{equation}
 respectively. The corresponding kernels are
\begin{equation}
 \begin{split}
  K_n(y',y)=&\oint\frac{dz'}{2\pi}\int_{-\infty}^\infty \frac{dz}{2\pi} \frac{1}{z'-z}\left(1-\left(\frac{z}{z'}\right)^n\right)e^{\gamma ({z'}^2-z^2)+\imath (y' z'-yz)}\prod_{j=1}^\infty e^{\imath\delta_j (z'-z)}\frac{1-\imath\delta_j z'}{1-\imath\delta_j z}
  \end{split}
\end{equation}
and
\begin{equation}
 \begin{split}
  K_n(y',y)=&\oint\frac{dz'}{2\pi}\int_{-\infty}^\infty \frac{dz}{2\pi} \frac{1}{z'-z}\left(1-\left(\frac{z}{z'}\right)^n\right)e^{\imath (y' z'-yz)}\prod_{j=1}^\infty \frac{1-\imath\delta_j z'}{1-\imath\delta_j z},
  \end{split}
\end{equation}
respectively. Hence these general results are relatively simple and explicit for these kinds of ensembles.

For the Gaussian unitary ensemble (GUE) we have $\omega(x)=\exp[-x^2/2]$. Then the bi-orthonormal functions become
\begin{equation}\label{gauss}
\begin{split}
p_j(y')=\oint\frac{dz'}{2\pi(\imath z')^{j+1}}e^{ {z'}^2/2+\imath y' z'}=\frac{1}{j!}H_j(y'),\ q_j(y)=\int_{-\infty}^\infty \frac{dz}{2\pi} (\imath z)^je^{- z^2/2-\imath yz}=\frac{1}{\sqrt{2\pi}}H_j(y)e^{-y^2/2}
\end{split}
\end{equation}
with $H_j$ the Hermite polynomials in the monic normalization.

Let us emphasize that P\'olya frequency functions of infinite order do not cover all P\'olya ensembles of a fixed matrix size. There are many more P\'olya frequency functions for example of the form
\begin{equation}\label{Laplace-Polya.c}
\begin{split}
\mathcal{F}\omega(s)=&\prod_{j=1}^\infty\frac{1}{(1-\imath \delta_j s)^{\nu_j}},\quad \nu_j,\delta_j\geq0,\, 0<\sum_{j=1}^\infty \nu_j\delta_j<\infty,\ {\rm and}\ {\rm Im}\, s>-\min\left\{\frac{1}{\delta_j}\right\},
 \end{split}
\end{equation}
and in the case that $\nu_j$ is not an integer it has to be $\nu_j>n-1$, see \cite{Kieburg:2017}. For $n=2$ we even only need to satisfy the condition that ${\rm log}\,\omega$ is concave which covers an extremely large class of functions.  The Laguerre ensemble obtained by the choice $\omega(x)=x^{n+\nu-1}\exp[-x]\Theta(x)$ with $\nu>-1$ and $\Theta(x)$ the Heaviside step function is of such a kind of P\'olya frequency function. Its Fourier transform is $\mathcal{F}\omega(s)=\Gamma[n+\nu](1-\imath s)^{-n-\nu}$. Hence the bi-orthonormal functions are
\begin{equation}\label{Lag}
\begin{split}
p_j(y')=&\frac{1}{\Gamma[n+\nu]}\oint\frac{dz'}{2\pi(\imath z')^{j+1}}e^{\imath y' z'}(1-\imath z')^{n+\nu}=\frac{1}{j!\Gamma[n+\nu]}L_j^{(n+\nu-j)}\left(y'\right),\\
 q_j(y)=&\Gamma[n+\nu]\int_{-\infty}^\infty \frac{dz}{2\pi} (\imath z)^je^{-\imath yz}(1-\imath z)^{-n-\nu}=L_j^{(n+\nu-j-1)}\left(y\right)y^{n+\nu-j-1}e^{-y}\Theta(y)
\end{split}
\end{equation}
with $L_k^{(\mu)}$ the generalized Laguerre polynomials in monic normalization. The reason why we obtain a different set of bi-orthonormal functions instead of the standard pair $\{L_j^{(\nu)}(y')/j!,L_j^{(\nu)}(y)y^\nu e^{-y}/\Gamma[j+\nu+1]\}_{j=0,\ldots,n-1}$, see~\cite{Mehta:2004,Forrester:2010,Akemann:2011} and end of Example~\ref{example.b}, follows from the fact that we looked for polynomials which are bi-orthonormal to the functions $q_j(y)=(-\partial_y)^{j} y^{n+\nu-1}e^{-y}\Theta(y)$.
\end{examples}

\subsection{Statistics with Fixed Matrices on $H_2$}\label{sec:evstat.fixed.a}

Again we state the result for general P\'olya ensembles and afterwards the results where we assume the same additional properties of the weight $\omega$ as in Corollary~\ref{cor:simpl.Polya.a}.

\begin{theorem}[Eigenvalue Statistics of P\'olya Ensembles with Fixed Matrices]\label{thm:stat.fixed.a}\

 The joint probability density~\eqref{jpdf.fixed.a} can be described by the bi-orthonormal functions
 \begin{equation}\label{biorth.fixed.a}
 \{p_j,q_j\}_{j=0,\ldots,n-1}=\left\{\int_{0}^\infty dr e^{-r}\left.\frac{(r-\imath\partial_t\partial_{y'})^{n-1}}{(n-1)!}\left(\prod_{l\neq j+1}\frac{x_l-y'}{x_l-x_{j+1}}\right)\frac{1}{\mathcal{F}\omega(t)}\right|_{t=0}\,,\,\omega(y-x_{j+1})\right\}_{j=0,\ldots,n-1}.
 \end{equation}
 The kernel is given by
 \begin{equation}\label{kernel.fixed.a}
 K_n(y',y)=\int_{0}^\infty dr e^{-r}\left.\frac{(r-\imath\partial_t\partial_{y'})^{n-1}}{(n-1)!}\left(\sum_{j=1}^n\omega(y-x_j)\prod_{l\neq j}\frac{x_l-y'}{x_l-x_j}\right)\frac{1}{\mathcal{F}\omega(t)}\right|_{t=0}.
 \end{equation}
\end{theorem}

Let us emphasize that this time the polynomials $p_j(y')$ are all of the same order namely of order $n-1$. The reason is the same as for the polynomials found in~\cite{Akemann:2015}. The weights $\omega(y-x_{j+1})$ only differ in the argument and thus the set of the weights is symmetric under permuting the variables $x$.

{\bf Proof:}\\
We show the bi-orthonormality by first noting that the polynomials of the Theorem~\ref{thm:stat.Polya.a}, which we denote now by  $\hat{p}_j(y')$, are given by
\begin{equation}\label{pol-proof.a}
\hat{p}_j(y')=\left.\frac{(y'-\imath\partial_t)^j}{j!}\frac{1}{\mathcal{F}\omega(t)}\right|_{t=0}=\int_{0}^\infty dr e^{-r}\left.\frac{(r-\imath\partial_t\partial_{y'})^{n-1}}{(n-1)!j!}\frac{y'^j}{\mathcal{F}\omega(t)}\right|_{t=0}.
\end{equation}
These polynomials are bi-orthonormal to $(-\partial_y)^j\omega(y)$.
The bi-orthonormality of $p_l(y')$ to $\omega(y-x_{m+1})$ is based on the bi-orthonormality of $\hat{p}_j(y')$ and $(-\partial_y)^j\omega(y)$  as can be shown as follows
\begin{equation}
\begin{split}
\int_{-\infty}^\infty dy p_l(y)\omega(y-x_{m})=&\int_{-\infty}^\infty dy\omega(y-x_{m})\int_{0}^\infty dr e^{-r}\left.\frac{(r-\imath\partial_t\partial_{y})^{n-1}}{(n-1)!}\left(\prod_{k\neq l}\frac{x_k-y}{x_k-x_l}\right)\frac{1}{\mathcal{F}\omega(t)}\right|_{t=0}\\
=&\int_{-\infty}^\infty dy\omega(y)\int_{0}^\infty dr e^{-r}\left.\frac{(r-\imath\partial_t\partial_{y})^{n-1}}{(n-1)!}\left(\prod_{k\neq l}\frac{x_k-x_m-y}{x_k-x_l}\right)\frac{1}{\mathcal{F}\omega(t)}\right|_{t=0}.
\end{split}
\end{equation}
In the next step we expand the product in monomials of $y$,
\begin{equation}
\begin{split}
\int_{0}^\infty dr e^{-r}\left.\frac{(r-\imath\partial_t\partial_{y})^{n-1}}{(n-1)!}\left(\prod_{k\neq l}\frac{x_k-x_m-y}{x_k-x_l}\right)\frac{1}{\mathcal{F}\omega(t)}\right|_{t=0}=&\sum_{k=0}^{n-1}c_k\int_{0}^\infty dr e^{-r}\left.\frac{(r-\imath\partial_t\partial_{y})^{n-1}}{(n-1)!}y^k\frac{1}{\mathcal{F}\omega(t)}\right|_{t=0}\\
=&\sum_{k=0}^{n-1}c_k\hat{p}_k(y).
\end{split}
\end{equation}
The coefficients $c_k$ are irrelevant apart from $c_0$ since the integral with $\omega(y)$ vanishes for all $\hat{p}_k(y)$ with $k>0$. The coefficient $c_0$ is $1$ for $l=m$  while it vanishes otherwise since then $x_m$ agrees with one of the $x_k$ in the product. Hence the bi-orthonormality follows and  the kernel is an  immediate consequence from the general definition~\eqref{k-point}.
\hfill$\square$
\vskip0.2cm

Now we come to a simplification similar to Corollary~\ref{cor:simpl.fixed.a}.

\begin{corollary}[Simplification of  Theorem~\ref{thm:stat.fixed.a}]\label{cor:simpl.fixed.a}\

We assume the requirements of  Theorem~\ref{thm:stat.fixed.a} and the properties of $\omega$ in Corollary~\ref{cor:simpl.Polya.a}. Then the bi-orthonormal functions can be written as
 \begin{equation}\label{biorth.fixed.a.2}
 \{p_j,q_j\}_{j=0,\ldots,n-1}=\left\{\oint\frac{dz'}{2\pi\imath {z'}}\int_0^\infty dx \frac{e^{\imath y' z'-x}}{\mathcal{F}\omega(z')}\left(\prod_{l\neq j+1}\frac{x_l+\imath x/z'}{x_l-x_{j+1}}\right)\,,\,\omega(y-x_{j+1})\right\}_{j=0,\ldots,n-1}
 \end{equation}
and the kernel as
 \begin{equation}\label{kernel.fixed.a.2}
 \begin{split}
  K_n(y',y)=&\oint\frac{dz'}{2\pi\imath {z'}}\int_0^\infty dx \frac{e^{\imath y' z'-x}}{\mathcal{F}\omega(z')}\left(\sum_{j=1}^n\omega(y-x_{j})\prod_{l\neq j}\frac{x_l+\imath x/z'}{x_l-x_{j}}\right).
  \end{split}
 \end{equation}
 The contour of $z'$ encircles the origin $z'=0$ counter clockwise.
\end{corollary}

{\bf Proof:}\\
Everything follows from the identity
\begin{equation}\label{proof.1.b}
\int_0^\infty dr e^{-r}\left.\frac{(r-\imath\partial_t\partial_{y'})^{n-1}}{(n-1)!}\frac{p(y')}{\mathcal{F}\omega(t)}\right|_{t=0}=\oint\frac{dz'}{2\pi\imath {z'}}\int_0^\infty dx \frac{e^{\imath y' z'-x}}{\mathcal{F}\omega(z')}p\left(\frac{x}{\imath z'}\right)
\end{equation}
for any polynomial $p(y')$ of a maximal order $n-1$. This can be readily checked by expanding the polynomial in monomials and evaluating the integrals. The integral over $x$ generates a factorial which would otherwise hinder the resummation to the exponential function $e^{\imath y' z'}$.
\hfill$\square$
\vskip0.2cm

We want to point out that similar products of the ratios in pairwise differences of $x_l$ and the integration variable was also found in similar settings, e.g. see~\cite{Forrester:2015,Akemann:2015b} where the sum was also expressed as contour integrals and~\cite{Guhr:1996a,Guhr:1996b} for the GUE derived via the supersymmetry method. Indeed when $\omega$ is holomorphic about the real line we can rewrite the sum in the bracket in Eq.~\eqref{kernel.fixed.a.2} as
\begin{equation}\label{contour:weight.a}
\sum_{j=1}^n\omega(y-x_{j})\prod_{l\neq j}\frac{x_l+\imath x/z'}{x_l-x_{j}}=-\oint\frac{dz}{2\pi\imath }\frac{\omega(y-z)}{z+\imath x/z'}\prod_{l=1}^n\frac{x_l+\imath x/z'}{x_l-z},
\end{equation}
where the contour only encircles the points $x_1,\ldots,x_n$ counter clockwise.
Note that the holomorphy of $\omega$ about the real axis is already not fulfilled for the Laguerre ensemble and one has to be careful with this formula while it is for the GUE, see the examples in~\ref{example.a.1}. Thus it is already a quite strong condition for P\'olya ensembles.

\begin{examples}\label{example.a.1}\
\begin{enumerate}
\item	Let us consider a random matrix $X_1$ drawn from the GUE, i.e. $\omega(x)=e^{-x^2/2}$. It is indeed well-known~\cite[Chapter 3]{Katori:2016} that due to the Harish-Chandra-Izykson-Zuber integral~\cite{Harish:1956,Itzykson:1980} the joint probability density function of the eigenvalues yields a polynomial ensemble on $D$ of the form~\eqref{jpdf.fixed.a}, especially that we have the weights $e^{-(y_a-x_b)^2/2}$ in the second determinant. The  formula~\ref{biorth.fixed.a.2} for the polynomials bi-orthonormal to these weights can be simplified by using the identity
 \begin{equation}\label{Gauss:id}
  \oint\frac{dz'}{2\pi\imath {z'}}\int_0^\infty \frac{dx}{\sqrt{2\pi}} e^{\imath y' z'-x+{z'}^2/2}\left(\frac{-\imath x}{z'}\right)^m=\int_{-\infty}^\infty\frac{dx'}{2\pi}e^{-{x'}^2/2}(y'+\imath x')^m
 \end{equation}
 valid for any $m=0,\ldots,n-1$. Then we have
 \begin{equation}
 p_j(y')=\int_{-\infty}^\infty\frac{dx'}{2\pi}e^{-{x'}^2/2}\left(\prod_{l\neq j+1}\frac{x_l-y'-\imath x'}{x_l-x_{j+1}}\right).
 \end{equation}
Also the kernel~\eqref{kernel.fixed.a.2} can be simplified for the particular case via the contour integral~\eqref{contour:weight.a} which yields
 \begin{equation}
  K_n(y',y)=-\int_{-\infty}^\infty\frac{dx'}{2\pi}\oint\frac{dz}{2\pi\imath }\frac{\exp[-(z^2+{x'}^2)/2]}{z-\imath x'}\prod_{l\neq j}\frac{x_l-y'-\imath x'}{x_l-y-z}.
 \end{equation}
 We substituted $z\to z+y$ in the contour integral in Eq.~\eqref{contour:weight.a} such that we enclose $z$ integral encloses the points $x_1+y,\ldots,x_n+y$ counter clockwise but it does not enclose $x'$. Also this result is well-known from the supersymmetry method when evaluating the regularization with an imaginary increment, see~\cite{Guhr:1996a,Guhr:1996b}.
\item	As a second example we want to consider the Laguerre ensemble with $\omega(x)=x^{n+\nu-1}e^{-x}\Theta(x)$. Again there is an identity namely
 \begin{equation}\label{Lag:id}
  \oint\frac{dz'}{2\pi\imath {z'}}\int_0^\infty dx e^{\imath y' z'-x}(1-\imath z')^{n+\nu}\left(\frac{-\imath x}{z'}\right)^m=\Gamma[n+\nu+1]\oint\frac{dx'}{2\pi\imath x'^{n+\nu+1}}e^{x'}(y'- x')^m,
 \end{equation}
 where $x'$ encircles the origin counter clockwise.
 Hence the polynomials become
 \begin{equation}
 p_j(y')=\Gamma[n+\nu+1]\oint\frac{dx'}{2\pi\imath x'^{n+\nu+1}}e^{x'}\left(\prod_{l\neq j+1}\frac{x_l-y'+x'}{x_l-x_{j+1}}\right).
 \end{equation}
 To use formula~\eqref{contour:weight.a} for the weight we have to assume that $y\neq x_l$ for any $l=1,\ldots,n$ since the weight is not holomorph at the origin but everywhere else it can be analytically continued in a neighbourhood about the real axis, in particular we can choose $\omega(z)=z^{n+\nu} e^{-z}\Theta({\rm Re}\, z)$ for $z\in\mathbb{C}$ with the imaginary axis as the non-homlomorphic set of this function. The kernel is then
 \begin{equation}
  K_n(y',y)=-\Gamma[n+\nu+1]\oint\frac{dx'}{2\pi \imath x'}\oint\frac{dz}{2\pi\imath }\Theta({\rm Re}\,z)\left(\frac{z}{x'}\right)^{n+\nu}\frac{\exp[-z+x']}{x'-z}\prod_{l\neq j}\frac{x_l-y'+x'}{x_l-y+z}
 \end{equation}
 for all $y\notin\{x_1,\ldots,x_n\}$, $x'$ encircles the origin counter clockwise and $z$ only encircles the points $x_1-y,\ldots,x_n-y$ (but not $z$) counter clockwise and closely enough such that we do not cross the non-holomorphic region. This time we substituted $z\to y-z$ in Eq.~\eqref{contour:weight.a} for the considered weight. Both results are completely new but have the flavor (on the structural level) that they can be derived by the supersymmetry method as well.
\end{enumerate}
\end{examples}

\subsection{Statistics with Polynomial Ensembles on $H_2$}\label{sec:evstat.poly.a}

Finally we let the formerly fixed matrix $X_2\in H_2$ be a random matrix, too. It shall be drawn from a polynomial ensemble. As before we state first the result for a random matrix $X_1\in H_2$ drawn from a general P\'olya ensemble on $H_2$.

\begin{theorem}[Eigenvalue Statistics of P\'olya Ensembles with Polynomial Ensembles]\label{thm:stat.poly.a}\

 Consider the joint probability density~\eqref{jpdf.poly.a}  where the polynomial ensemble of $X_2$ can be described by the bi-orthonormal functions $\{\tilde{p}_j,w_{j+1}\}_{j=0,\ldots,n-1}$ and its kernel is $\tilde{K}_n(y'y)=\sum_{j=0}^{n-1}\tilde{p}_j(y')w_{j+1}(y)$. The pair of bi-orthonormal functions corresponding to $Y=X_1+X_2$ is
 \begin{equation}\label{biorth.poly.a}
 \{p_j,q_j\}_{j=0,\ldots,n-1}=\left\{\int_{0}^\infty dr e^{-r}\left.\frac{(r-\imath\partial_t\partial_{y'})^{n-1}}{(n-1)!}\frac{\tilde{p}_j(y')}{\mathcal{F}\omega(t)}\right|_{t=0}\,,\,\omega\ast w_{j+1}(y)\right\}_{j=0,\ldots,n-1}
 \end{equation}
and the corresponding kernel is
 \begin{equation}\label{kernel.poly.a}
 \begin{split}
 K_n(y',y)=&\int_{-\infty}^\infty d\hat{y}\int_{0}^\infty dr e^{-r}\omega(y-\hat{y})\left.\frac{(r-\imath\partial_t\partial_{y'})^{n-1}}{(n-1)!}\frac{\tilde{K}_n(y',\hat{y})}{\mathcal{F}\omega(t)}\right|_{t=0}.
 \end{split}
 \end{equation}
\end{theorem}

{\bf Proof:}\\
The proof works along the same ideas as the one of Theorem~\ref{thm:stat.fixed.a}. The bi-orthonormality of $ \{p_j,q_j\}_{j=0,\ldots,n-1}$ is again based on the bi-orthonormality of $ \{\hat{p}_j(y'),(-\partial_y)^j\omega(y)\}_{j=0,\ldots,n-1}$ with $\hat{p}_j$ as in Eq.~\eqref{pol-proof.a}. The integral we have to consider is
\begin{equation}
\begin{split}
\int_{-\infty}^\infty dy p_l(y)\omega\ast w_{m+1}(y)=&\int_{-\infty}^\infty dy\int_{-\infty}^\infty dx\omega(y-x)w_{m+1}(x)\int_{0}^\infty dr e^{-r}\left.\frac{(r-\imath\partial_t\partial_{y})^{n-1}}{(n-1)!}\frac{\tilde{p}_l(y)}{\mathcal{F}\omega(t)}\right|_{t=0}\\
=&\int_{-\infty}^\infty dy\int_{-\infty}^\infty dx\omega(y)w_{m+1}(x)\int_{0}^\infty dr e^{-r}\left.\frac{(r-\imath\partial_t\partial_{y})^{n-1}}{(n-1)!}\frac{\tilde{p}_l(y+x)}{\mathcal{F}\omega(t)}\right|_{t=0}.
\end{split}
\end{equation}
We underline that everything is absolutely integrable such that we can interchange the integrals without any problems.
This time we expand $\tilde{p}_l(y+x)=\sum_{k=0}^ld_k(x) y^k$ with $d_k(x)$ some polynomials in $x$. The integral over $r$ together with the derivatives yields the polynomials $\hat{p}_k(y)$. After integrating over $y$ with the weight $\omega(y)$ only the term $d_{k=0}(x)=\tilde{p}_l(x)$ survives which is bi-orthonormal to $w_{m+1}(y)$ which proves the bi-orthonormality of $ \{p_j,q_j\}_{j=0,\ldots,n-1}$. The kernel is again a direct result because the finite sum in Eq.~\eqref{k-point} can be interchanged with the integrals and derivatives.
\hfill$\square$
\vskip0.2cm

Also for this case a simplification exists when $\omega$ satisfies the additional conditions of Corollary~\ref{cor:simpl.Polya.a}.

\begin{corollary}[Simplification of  Theorem~\ref{thm:stat.poly.a}]\label{cor:simpl.poly.a}\

We assume the requirements of  Theorem~\ref{thm:stat.poly.a} and the properties of $\omega$ in Corollary~\ref{cor:simpl.Polya.a}. Then the bi-orthonormal functions can be written as
 \begin{equation}\label{biorth.poly.a.2}
 \{p_j,q_j\}_{j=0,\ldots,n-1}=\left\{\oint\frac{dz'}{2\pi\imath {z'}}\int_0^\infty dx \frac{e^{\imath y' z'-x}}{\mathcal{F}\omega(z')}\tilde{p}_j\left(\frac{x}{\imath z'}\right)\,,\,\omega\ast w_{j+1}(y)\right\}_{j=0,\ldots,n-1}
 \end{equation}
and the kernel as
 \begin{equation}\label{kernel.poly.a.2}
  K_n(y',y)=\int_{-\infty}^\infty d\hat{y}\oint\frac{dz'}{2\pi\imath {z'}}\int_0^\infty dx \frac{e^{\imath y' z'-x}}{\mathcal{F}\omega(z')}\omega(y-\hat{y})\tilde{K}_n\left(\frac{x}{\imath z'},\hat{y}\right).
 \end{equation}
 The contour of $z'$ encircles the origin $z'=0$ counter clockwise.
\end{corollary}

{\bf Proof:}\\
Again relation~\eqref{proof.1.b} proves these statements.
\hfill$\square$
\vskip0.2cm

Unfortunately, we were not able to reduce the kernel~\eqref{kernel.poly.a.2} to  an expression with less integrals to perform for a general P\'olya ensemble on $H_2$. Certainly for particular cases like the convolution with a GUE such a simplification exists as shown in~\cite{Claeys:2015}, see also the following examples.

\begin{examples}\label{example.a.2}\
\begin{enumerate}
 \item	Let us again start with the GUE. We use identity~\eqref{Gauss:id} and find for the polynomials
 \begin{equation}
  p_j(y')=\int_{-\infty}^\infty\frac{dx'}{2\pi}e^{-{x'}^2/2}\tilde{p}_j\left(y'+\imath x'\right)=\int_{-\infty}^\infty\frac{dx'}{\sqrt{2\pi}}e^{(\imath x'- y')^2/2}\tilde{p}_j\left(\imath x'\right).
 \end{equation}
 Then the kernel~\eqref{kernel.poly.a.2} becomes in the Gaussian case
 \begin{equation}
  K_n(y',y)=\int_{-\infty}^\infty d\hat{y}\int_{-\infty}^\infty\frac{dx'}{2\pi}\exp\left[-\frac{(\hat{y}-y)^2}{2}+\frac{(\imath x'- y')^2}{2}\right]\tilde{K}_n\left(\imath x',\hat{y}\right).
 \end{equation}
 This result was already derived in~\cite{Claeys:2015}. This transformation formula was applied in~\cite{Claeys:2016} where the polynomial ensemble was chosen to be the Laguerre ensemble, products of Ginibre matrices or matrices drawn from the Jacobi ensemble and Muttalib-Borodin ensembles.
 \item	As a second example we again consider the Laguerre ensemble and make use of the formula~\eqref{Lag:id}. The polynomials are
 \begin{equation}
 p_j(y')=\Gamma[n+\nu+1]\oint\frac{dx'}{2\pi\imath x'^{n+\nu+1}}e^{x'}\tilde{p}_j\left(y'-x'\right)
 \end{equation}
 and the kernel is
 \begin{equation}
  K_n(y',y)=\Gamma[n+\nu+1]\int_{0}^\infty \frac{d\hat{y}}{\hat{y}}\oint\frac{dx'}{2\pi \imath x'}\left(\frac{\hat{y}}{x'}\right)^{n+\nu}\exp\left[-\hat{y}+x'\right]\tilde{K}_n\left(y'-x',y-\hat{y}\right)
 \end{equation}
 as can be easily checked. These results were already found in~\cite{Kuijlaars:2016a}.
\end{enumerate}
\end{examples}

A mixture of the two examples with the results in subsection~\ref{sec:evstat.fixed.a} was considered in~\cite[Appendix E]{Forrester:2005}. There the random matrix $A+XBX^*+\sqrt{t} Y$ was studied with $t\in\mathbb{R}_+$, $A,B\in H_2$ fixed, a Hermitian random matrix $Y\in H_2$ drawn from a GUE and  a complex random matrix $X\in{\rm Gl}_{\mathbb{C}}(n)$ drawn from a Laguerre (Wishart) ensemble. For $B=0$ we have the shift of a GUE with a fixed matrix, namely $A$. For $A=0$, we have the situation of the GUE shifted by the polynomial ensemble described by  the random matrix $XBX^*$.

\section{Convolutions with P\'olya Ensembles on $M=H_1,H_4,M_\nu$}\label{sec:M}

We follow the same lines as in Sec.~\ref{sec:H2} by first deriving the bi-orthonormal functions of a P\'olya ensemble on $M=H_1,H_4,M_\nu$ without any shift in subsection~\ref{sec:evstat.Polya.b}, and then stating and proving the results for a shift with a fixed matrix in $M$ (subsection~\ref{sec:evstat.fixed.b}) and with a polynomial ensemble on $M$ (subsection~\ref{sec:evstat.poly.b}). We also consider examples as we have done for the additive convolution on Hermitian matrices. These examples can be found in~\ref{example.b}.

\subsection{Statistics of P\'olya Ensembles on $M=H_1,H_4,M_\nu$}\label{sec:evstat.Polya.b}

As before we start with the general case where we look for the bi-orthonormal functions of a random matrix in $M=H_1,H_4,M_\nu$ drawn from a P\'olya ensemble on $M$ but completely without a shift.

\begin{theorem}[Eigenvalue Statistics of P\'olya Ensembles]\label{thm:stat.Polya.b}\

 The joint probability density~\eqref{jpdf.Polya.b} can be described by the bi-orthonormal functions
 \begin{equation}\label{biorth.Polya.b}
 \{p_j,q_j\}_{j=0,\ldots,n-1}=\left\{\int_0^\infty dr e^{-r}\left.\frac{\Gamma(\nu+1)(r-\partial_t{y'}^{-\nu}\partial_{y'}{y'}^{\nu+1}\partial_{y'})^{n-1}}{(n-1)!j!\Gamma(\nu+j+1)}\frac{{y'}^j}{\mathcal{H}_\nu\omega(t)}\right|_{t=0}\,,\,(\partial_yy^{\nu+1}\partial_yy^{-\nu})^j\omega(y)\right\}_{j=0,\ldots,n-1}.
 \end{equation}
  The kernel is given by
 \begin{equation}\label{kernel.Polya.b}
 \begin{split}
  K_n(y',y)=&\iint_0^\infty drdR R^{n-1} e^{-r-R}\left.\frac{\Gamma(\nu+1)(r-\partial_t{y'}^{-\nu}\partial_{y'}{y'}^{\nu+1}\partial_{y'})^{n-1}}{(n-1)!j!\Gamma(\nu+j+1)}\frac{1}{\mathcal{H}_\nu\omega(t)}\right|_{t=0}L_{n-1}^{(\nu)}\left(\frac{y'}{R}\partial_yy^{\nu+1}\partial_yy^{-\nu}\right)\omega(y)
  \end{split}
 \end{equation}
 with $L_j^{(\nu)}(x)=x^j+\ldots$ the generalized Laguerre polynomials in monic normalization.
\end{theorem}

{\bf Proof:}\\
Let us emphasize that all polynomials of order $m<l$ are automatically orthogonal to $(\partial_yy^{\nu+1}\partial_yy^{-\nu})^m\omega(y)$ as can be readily checked via integration by parts. Thus we have to construct the polynomials in such a way that they become orthogonal also to the weights with $m>l$. The bi-orthonormal polynomials to the weights $\{(\partial_yy^{\nu+1}\partial_yy^{-\nu})^j\omega(y)\}_{j=0,\ldots,n-1}$ can be constructed via taking the following determinants,
\begin{equation}
p_m(y)=\det\left[\begin{array}{c} \int_{0}^\infty \hat{y}^{b-1} (\partial_{\hat{y}}\hat{y}^{\nu+1}\partial_{\hat{y}}\hat{y}^{-\nu})^{a-1}\omega(\hat{y}) d\hat{y} \\ y^{b-1} \end{array}\right]_{\substack{a=1,\ldots,m\\b=1,\ldots,m+1}}\biggl/\det\left[\begin{array}{c} \int_{0}^\infty \hat{y}^{b-1} (\partial_{\hat{y}}\hat{y}^{\nu+1}\partial_{\hat{y}}\hat{y}^{-\nu})^{a-1}\omega(\hat{y}) d\hat{y}\end{array}\right]_{a,b=1,\ldots,m+1}.
\end{equation}
The one-dimensional integrals are
\begin{equation}
\begin{split}
\int_{0}^\infty \hat{y}^{b-1} (\partial_{\hat{y}}\hat{y}^{\nu+1}\partial_{\hat{y}}\hat{y}^{-\nu})^{a-1}\omega(\hat{y}) d\hat{y}=&\Theta_{ab}\frac{(b-1)!\Gamma(\nu+b)}{(b-a)!\Gamma(\nu+b-a+1)}\int_{0}^\infty \hat{y}^{b-a}\omega(\hat{y}) d\hat{y}\\
=&\Theta_{ab}\frac{(b-1)!\Gamma(\nu+b)}{(b-a)!\Gamma(\nu+1)}\left.(-\partial)_t^{b-a}\mathcal{H}_\nu\omega(t)\right|_{t=0}
\end{split}
\end{equation}
with $\Theta_{ab}$ as before, namely $\Theta_{ab}=1$ for $a\leq b$ and $\Theta_{ab}=0$ for $a>b$. In the second step we employed the identity
\begin{equation}
\int_{0}^\infty \hat{y}^{k}\omega(\hat{y}) d\hat{y}=\frac{\Gamma(\nu+k+1)}{\Gamma(\nu+1)}\left.(-\partial)_t^{k}\mathcal{H}_\nu\omega(t)\right|_{t=0},\quad{for\ any\ }k\in\mathbb{N}_0,
\end{equation}
because of the integrability conditions of $\omega\in L_\nu^1(\mathbb{R}_+)$, see Eq.~\eqref{L1-sets}. Hence we have
\begin{equation}
p_m(y)=\frac{1}{\mathcal{H}_\nu\omega(0)}\det\left[\begin{array}{c} \displaystyle\Theta_{ab}\left.\frac{(-\partial_t)^{b-a}}{(b-a)!}\frac{\mathcal{H}_\nu\omega(t)}{\mathcal{H}_\nu\omega(0)}\right|_{t=0} \\ \displaystyle\frac{\Gamma(\nu+1)}{(b-1)!\Gamma(\nu+b)}y^{b-1} \end{array}\right]_{\substack{a=1,\ldots,m\\b=1,\ldots,m+1}}.
\end{equation}
When expanded this expression becomes
\begin{equation}\label{proof1.a}
p_m(y)=\frac{(-1)^m}{\mathcal{H}_\nu\omega(0)}\sum_{j=0}^m\frac{\Gamma(\nu+1)(-y)^j}{j!\Gamma(\nu+j+1)}\det\left[\begin{array}{c} \displaystyle\Theta_{ab}\left.\frac{(-\partial_t)^{b-a}}{(b-a)!}\frac{\mathcal{H}_\nu\omega(t)}{\mathcal{H}_\nu\omega(0)}\right|_{t=0} \end{array}\right]_{\substack{a=j+1,\ldots,m\\ b=j+2,\ldots,m+1}}.
\end{equation}
The remaining determinants are a Toeplitz determinants of the form~\eqref{toeplitz} with $L=1$, $n=m-j$ and $c_{j}=1/(j+1)!(-\partial_t)^{j+1}\mathcal{H}_\nu\omega(t)/\mathcal{H}_\nu\omega(0)|_{t=0}$. Thus we have
\begin{equation}
\det\left[\begin{array}{c} \displaystyle\Theta_{ab}\left.\frac{(-\partial_t)^{b-a}}{(b-a)!}\frac{\mathcal{H}_\nu\omega(t)}{\mathcal{H}_\nu\omega(0)}\right|_{t=0} \end{array}\right]_{\substack{a=j+1,\ldots,m\\b=j+2,\ldots,m+1}}=\left.\frac{\partial_t^{m-j}}{(m-j)!}\frac{\mathcal{H}_\nu\omega(0)}{\mathcal{H}_\nu\omega(t)}\right|_{t=0}
\end{equation}
which we plug into the sum~\eqref{proof1.a}. It can be readily checked that this sum can be rewritten to
\begin{equation}
p_m(y)=\int_0^\infty dr e^{-r}\left.\frac{\Gamma(\nu+1)(r-\partial_t{y}^{-\nu}\partial_{y}{y}^{\nu+1}\partial_{y})^{n-1}}{(n-1)!m!\Gamma(\nu+m+1)}\frac{{y}^m}{\mathcal{H}_\nu\omega(t)}\right|_{t=0}.
\end{equation}
The sum for the kernel is equal to the generalized Laguerre polynomial up to a factor $1/(n-1-j)!$. This factor is introduced by the integral over $R$ in Eq.~\eqref{kernel.Polya.b}.\hfill$\square$
\vskip0.2cm

As for Theorem~\ref{thm:stat.Polya.a} we can simplify the results when we assume some additional conditions for the weight $\omega$.

\begin{corollary}[Simplification of  Theorem~\ref{thm:stat.Polya.b}]\label{cor:simpl.Polya.b}
Additionally to the requirements of  Theorem~\ref{thm:stat.Polya.b} we assume that the Hankel transform $\mathcal{H}_\nu\omega$ is holomorphic at the origin and $z^{j+(2\nu-1)/4}\mathcal{H}_\nu\omega(z)$ is absolutely integrable along the positive real line for $j=0,\ldots,n-1$. Then the bi-orthonormal functions simplify to
 \begin{equation}\label{biorth.Polya.b.2}
 \{p_j,q_j\}_{j=0,\ldots,n-1}=\left\{\oint\frac{dz'}{2\pi\imath {z'}^{j+1}}\frac{I_\nu(2\sqrt{y'z'})}{(y'z')^{\nu/2}}\frac{\Gamma(\nu+1)}{\mathcal{H}_\nu\omega(z')}\,,\,\int_0^\infty dzz^jJ_\nu(2\sqrt{yz})(yz)^{\nu/2}\frac{\mathcal{H}_\nu\omega(z)}{\Gamma(\nu+1)}\right\}_{j=0,\ldots,n-1}.
 \end{equation}
and the kernel simplifies to
 \begin{equation}\label{kernel.Polya.b.2}
 \begin{split}
  K_n(y',y)=&\oint\frac{dz'}{2\pi\imath}\int_0^\infty dz\frac{1}{z'-z}\left(1-\left(\frac{z}{z'}\right)^n\right)\left(\frac{yz}{y'z'}\right)^{\nu/2}I_\nu(2\sqrt{y'z'})J_\nu(2\sqrt{yz})\frac{\mathcal{H}_\nu\omega(z)}{\mathcal{H}_\nu\omega(z')}.
  \end{split}
 \end{equation}
 The contour of $z'$ encircles the origin $z'=0$ counter clockwise while the contour for $z$ is along the positive real line.
\end{corollary}

{\bf Proof:}\\
As before for the proof of Corollary~\ref{cor:simpl.Polya.a} the absolute integrability of  $z^{j+(2\nu-1)/4}\mathcal{H}_\nu\omega(z)$ for $j=0,\ldots,n-1$ and the holomorphy of $\mathcal{H}_\nu\omega(z')$ at $z'=0$ with $\mathcal{H}_\nu\omega(0)\neq0$ guarantee us to rewrite the weights $q_j$ and the polynomials into the integral forms~\eqref{biorth.Polya.b.2}. For this purpose we emphasize that we first write the polynomials as a sum where the integral over $r$ in Eq.~\eqref{biorth.Polya.b} and the derivatives in $y'$ are evaluated. Then it is very simple to see that the contour integral in Eq.~\eqref{biorth.Polya.b.2} generate the same coefficients.  The sum for the kernel is the geometric sum and readily yields~\eqref{kernel.Polya.b.2}.
\hfill$\square$
\vskip0.2cm

Again we want to underline the resemblance to supersymmetry results though without Bessel functions, cf.~\cite{Guhr:2006,Kieburg:2009}. Even the Bessel function could be anticipated as they are known from the supersymmetric Berezin-Karpelevic integral~\cite{Guhr:1996} which is directly related to the complex Laguerre ensemble, see~\cite{Jackson:1996}.

\begin{examples}\label{example.b}
As a conclusion of this subsection we want to give a quite general example for the P\'olya ensembles on $M=H_1,H_4,M_\nu$. We use the sufficient condition for this kind of P\'olya ensembles derived in~\cite{Kieburg:2017} which relates the P\'olya ensembles on $M$ with those on $H_2$ where the support is only on the positive definite matrices. Let us assume that $\tilde\omega\in L^1_{\mathcal{F}}(\mathbb{R})$ corresponds to a P\'olya ensemble on $H_2$ for positive definite matrices, i.e. $\tilde\omega$ has only support on $\mathbb{R}_+$. Then the weight
\begin{equation}\label{Polya.rel}
\omega(x)=\frac{1}{\Gamma[\nu+1]}\int_0^\infty \left(\frac{x}{y}\right)^\nu \exp\left[-\frac{x}{y}\right]\tilde\omega(y)\frac{dy}{y}\in L^{1}_{\nu}(\mathbb{R}_+)
\end{equation}
belongs to a P\'olya ensemble on $M$. The Hankel transform of $\omega(x)$ is then
\begin{equation}
\mathcal{H}_\nu\omega(s)=\int_0^\infty \exp\left[- ys\right]\tilde\omega(y)dy.
\end{equation}
Indeed the integral for the Hankel transform and the integral over $y$ can be interchanged since they are absolutely integrable. Thus the Hankel transform of $\omega$ is equal to the Laplace transform of $\tilde\omega$. Choosing $\tilde\omega$ as a P\'olya frequency function of infinite order it can have only the form
\begin{equation}\label{Laplace-Polya.d}
\begin{split}
\mathcal{H}_\nu\omega(s)=&\mathcal{F}\tilde\omega(\imath s)=e^{-\delta s}\prod_{j=1}^\infty\frac{1}{1+\delta_j s},\quad \delta,\delta_j\geq0,\, 0<\sum_{j=1}^\infty \delta_j<\infty,\ {\rm and}\ {\rm Im}\, s>-\min\left\{\frac{1}{\delta_j}\right\},
 \end{split}
\end{equation}
cf. Eq.~\eqref{Laplace-Polya.b}. This time we cannot shift the factor $e^{-\delta s}$ since origin is an exceptional point in the spectrum which cannot be crossed. To satisfy the correct differentiability conditions we need that at least $n+1$ of the $\delta_j$ are non-vanishing. 

The bi-orthonormal pair of functions are for this kind of P\'olya ensembles
 \begin{equation}
 \begin{split}
 \{p_j,q_j\}_{j=0,\ldots,n-1}=&\biggl\{\Gamma(\nu+1)\oint\frac{dz'}{2\pi\imath {z'}^{j+1}}\frac{I_\nu(2\sqrt{y'z'})}{(y'z')^{\nu/2}}e^{\delta z'}\prod_{l=1}^\infty(1+\delta_l z')\,,\\
 &\frac{1}{\Gamma(\nu+1)}\int_0^\infty dzz^jJ_\nu(2\sqrt{yz})(yz)^{\nu/2}e^{-\delta z}\prod_{l=1}^\infty\frac{1}{1+\delta_l z}\biggl\}_{j=0,\ldots,n-1}.
 \end{split}
 \end{equation}
and the kernel becomes
 \begin{equation}
 \begin{split}
  K_n(y',y)=&\oint\frac{dz'}{2\pi\imath}\int_0^\infty dz\frac{1}{z'-z}\left(1-\left(\frac{z}{z'}\right)^n\right)\left(\frac{yz}{y'z'}\right)^{\nu/2}I_\nu(2\sqrt{y'z'})J_\nu(2\sqrt{yz})e^{\delta (z'-z)}\prod_{j=1}^\infty\frac{1+\delta_j z'}{1+\delta_j z}.
  \end{split}
 \end{equation}
 We want to underline that P\'olya frequency functions with a Fourier transform of the form~\eqref{Laplace-Polya.a} do not correspond to P\'olya ensembles on $M=H_1,H_4,M_\nu$ because they have a support on whole $\mathbb{R}$. However we can choose P\'olya frequency functions of finite order $n$ as given in Eq.~\eqref{Laplace-Polya.c} yielding weights $\omega$ via the relation~\eqref{Polya.rel} which satisfy all the required conditions of Corollary~\eqref{cor:simpl.Polya.b}.
 
 Finally we want to point out that not all P\'olya ensembles on $M$ can be found by Eq.~\eqref{Polya.rel} by choosing $\tilde{\omega}$ a P\'olya frequency function. For example choosing $\tilde{\omega}(y)=\Gamma[\nu+1]\delta(y-1)$ a Dirac delta function we recover the Laguerre ensemble now represented by the weight $\omega(x)=x^\nu e^{-x}$. We want to underline that the Laguerre ensemble discussed in Examples~\ref{example.a},~\ref{example.a.1}, and \ref{example.a.2} is the same as here. Only the weight is different since the differential operators we apply in the determinants~\eqref{jpdf.Polya.a} and~\eqref{jpdf.Polya.b} are different. Its Hankel transform is simply $\mathcal{H}_\nu\omega(s)=\Gamma[\nu+1]e^{-s}$. The polynomials and weights are
 \begin{equation}
 \begin{split}
 p_j(y')=&\oint\frac{dz'}{2\pi\imath {z'}^{j+1}}\frac{I_\nu(2\sqrt{y'z'})}{(y'z')^{\nu/2}}e^{z'}=\frac{1}{j!\Gamma[\nu+j+1]}L_j^{(\nu)}(y'),\\
 q_j(y)=&\int_0^\infty dzz^jJ_\nu(2\sqrt{yz})(yz)^{\nu/2}e^{-z}=L_j^{(\nu)}(y)y^\nu e^{-y}.
 \end{split}
 \end{equation}
 and are thus the standard choice of the Laguerre ensembles as can be found everywhere in the literature~\cite{Mehta:2004,Forrester:2010,Akemann:2011}.
\end{examples}

\subsection{Statistics with Fixed Matrices on $M=H_1,H_4,M_\nu$}\label{sec:evstat.fixed.b}

The theorem for the shifted P\'olya Ensemble on $M=H_1,H_4,M_\nu$ looks slightly more complicated than without the shift since do not get a unified formula because of the case $\nu=-1/2$. Nonetheless the ideas work a long the same lines as before.

\begin{theorem}[Eigenvalue Statistics of P\'olya Ensembles with Fixed Matrices]\label{thm:stat.fixed.b}\

 The joint probability density~\eqref{jpdf.fixed.b} can be described by the bi-orthonormal functions
 \begin{equation}\label{biorth.fixed.b}
 \begin{split}
 p_j(y')=&\int_{0}^\infty dr e^{-r}\left.\frac{(r-\partial_s{y'}^{-\nu}\partial_{y'}{y'}^{\nu+1}\partial_{y'})^{n-1}}{(n-1)!}\left(\prod_{l\neq j+1}\frac{x_l-y'}{x_l-x_{j+1}}\right)\frac{1}{\mathcal{H}_\nu\omega(s)}\right|_{s=0},\\
 q_j(y)=&\left\{\begin{array}{cl} \displaystyle\frac{\Gamma[\nu+1]}{\sqrt{\pi}\Gamma[\nu+1/2]}y^\nu\int_{-1}^{1} \frac{\omega(y+x_{j+1}-2\sqrt{yx_{j+1}}t)}{(y+x_{j+1}-2\sqrt{yx_{j+1}}t)^{\nu}}(1-t^2)^{\nu-1/2}dt, & {\rm for}\ \nu>-1/2, \\ \displaystyle \frac{1}{2\sqrt{y}}\left[ |\sqrt{y}-\sqrt{x_{j+1}}|\omega((\sqrt{y}-\sqrt{x_{j+1}})^2)+|\sqrt{y}+\sqrt{x_{j+1}}|\omega((\sqrt{y}+\sqrt{x_{j+1}})^2)\right], & {\rm for}\ \nu=-1/2 \end{array}\right.
 \end{split}
 \end{equation}
 with $j=0,\ldots,n-1$.
 The kernel is given by
 \begin{equation}\label{kernel.fixed.b}
 \begin{split}
 K_n(y',y)=&y^\nu\int_{0}^\infty dr e^{-r}\frac{(r-\partial_s{y'}^{-\nu}\partial_{y'}{y'}^{\nu+1}\partial_{y'})^{n-1}}{(n-1)!}\frac{1}{\mathcal{H}_\nu\omega(s)}\biggl|_{s=0}\\
 &\times\left(\sum_{j=1}^n\frac{\Gamma[\nu+1]}{\sqrt{\pi}\Gamma[\nu+1/2]}\int_{-1}^{1} \frac{\omega(y+x_{j}-2\sqrt{yx_{j}}t)}{(y+x_{j}-2\sqrt{yx_{j}}t)^{\nu}}(1-t^2)^{\nu-1/2}dt\prod_{l\neq j}\frac{x_l-y'}{x_l-x_j}\right)
 \end{split}
 \end{equation}
 for $\nu>-1/2$ and
 \begin{equation}\label{kernel.fixed.c}
 \begin{split}
 K_n(y',y)=&\frac{1}{2\sqrt{y}}\int_{0}^\infty dr e^{-r}\frac{(r-\partial_s{y'}^{-\nu}\partial_{y'}{y'}^{\nu+1}\partial_{y'})^{n-1}}{(n-1)!}\frac{1}{\mathcal{H}_\nu\omega(s)}\biggl|_{s=0}\\
 &\times\left(\sum_{j=1}^n\left[ |\sqrt{y}-\sqrt{x_{j}}|\omega((\sqrt{y}-\sqrt{x_{j}})^2)+|\sqrt{y}+\sqrt{x_{j}}|\omega((\sqrt{y}+\sqrt{x_{j}})^2)\right]\prod_{l\neq j}\frac{x_l-y'}{x_l-x_j}\right)
 \end{split}
 \end{equation}
 for $\nu=-1/2$.
\end{theorem}

{\bf Proof:}\\
Also this theorem is proved by explicitly showing the bi-orthogonality, i.e. that
\begin{equation}
\begin{split}
I_{ml}=&\frac{\Gamma[\nu+1]}{\sqrt{\pi}\Gamma[\nu+1/2]}\int_0^\infty dy p_m(y)y^\nu\left(\int_{-1}^{1} \frac{\omega(y+x_{l+1}-2\sqrt{yx_{l+1}}t)}{(y+x_{l+1}-2\sqrt{yx_{l+1}}t)^{\nu}}(1-t^2)^{\nu-1/2}dt\right)\\
=&\int_{\mathbb{C}^{\nu+1}} dv p_m(||v||^2)\frac{\omega(||v-\sqrt{x_{l+1}}e_1||^2)}{||v-\sqrt{x_{l+1}}e_1||^{2\nu}}\\
=&\int_{\mathbb{C}^{\nu+1}} dv p_m(||v+\sqrt{x_{l+1}}e_1||^2)\frac{\omega(||v||^2)}{||v||^{2\nu}}\\
=&\frac{\Gamma[\nu+1]}{\sqrt{\pi}\Gamma[\nu+1/2]}\int_0^\infty dy \left(\int_{-1}^1p_m(y+x_{l+1}+2\sqrt{yx_{l+1}}t)(1-t^2)^{\nu-1/2}dt\right)\omega(y)
\end{split}
\end{equation}
is equal to the Kronecker symbol. The equation above is only true for $\nu\in\{+1/2\}\cup\mathbb{N}_0$. For $\nu=-1/2$ we have
\begin{equation}
\begin{split}
I_{ml}=&\frac{1}{2}\int_0^\infty \frac{dy}{\sqrt{y}} p_m(y)\left[ |\sqrt{y}-\sqrt{x_{l+1}}|\omega((\sqrt{y}-\sqrt{x_{l+1}})^2)+|\sqrt{y}+\sqrt{x_{l+1}}|\omega((\sqrt{y}+\sqrt{x_{l+1}})^2)\right]\\
=&\frac{1}{2}\int_0^\infty dy\left[ p_m((\sqrt{y}-\sqrt{x_{l+1}})^2)+p_m((\sqrt{y}+\sqrt{x_{l+1}})^2)\right]\omega(y).
\end{split}
\end{equation}
We denote by $||v||$ the Euclidean norm and by $e_1$ a unit vector in $\mathbb{C}^{\nu+1}$.

In the next step we will show that the differential operator ${y}^{-\nu}\partial_{y}{y}^{\nu+1}\partial_{y}$ interchanges with the integral over $t$ for $\nu>-1/2$. Since $p_j$ is a polynomial it is enough to show it for monomials, in particular we have to show
\begin{equation}\label{proof2.a}
\int_{-1}^1m(m+\nu)(y+x_{l+1}+2\sqrt{yx_{l+1}}t)^{m-1}(1-t^2)^{\nu-1/2}dt={y}^{-\nu}\partial_{y}{y}^{\nu+1}\partial_{y}\int_{-1}^1(y+x_{l+1}+2\sqrt{yx_{l+1}}t)^{m}(1-t^2)^{\nu-1/2}dt.
\end{equation}
The left hand side is exactly the operator ${y}^{-\nu}\partial_{y}{y}^{\nu+1}\partial_{y}$ applied to $y^m$ and then replacing $y\to y+x_{l+1}+2\sqrt{yx_{l+1}}t$ and integrating over $t$. To show this we rewrite the integral as follows
\begin{equation}
\frac{1}{\sqrt{\pi}\Gamma[\nu+1/2]}\int_{-1}^1(y+x_{l+1}+2\sqrt{yx_{l+1}}t)^{m}(1-t^2)^{\nu-1/2}dt=m!\oint\frac{dz}{2\pi\imath z^{m+1}}e^{(y+x_{l+1})z}\frac{I_\nu(2\sqrt{yx_{l+1}}z)}{(\sqrt{yx_{l+1}}z)^\nu},
\end{equation}
where $I_\nu$ is the modified Bessel function of the first case and the contour encircles the origin. The differential operator can be evaluated as follows
\begin{equation}
\begin{split}
&{y}^{-\nu}\partial_{y}{y}^{\nu+1}\partial_{y}\oint\frac{dz}{2\pi\imath z^{m+1}}e^{(y+x_{l+1})z}\frac{I_\nu(2\sqrt{yx_{l+1}}z)}{(\sqrt{yx_{l+1}}z)^\nu}\\
=&\oint\frac{dz}{2\pi\imath z^{m+1}}e^{(y+x_{l+1})z}\left(yz^2+2zy\partial_y+(\nu+1)z+{y}^{-\nu}\partial_{y}{y}^{\nu+1}\partial_{y}\right)\frac{I_\nu(2\sqrt{yx_{l+1}}z)}{(\sqrt{yx_{l+1}}z)^\nu}\\
=&\oint\frac{dz}{2\pi\imath z^{m+1}}e^{(y+x_{l+1})z}\left(yz^2+2zy\partial_y+(\nu+1)z+x_{l+1}z^2\right)\frac{I_\nu(2\sqrt{yx_{l+1}}z)}{(\sqrt{yx_{l+1}}z)^\nu}\\
=&\oint\frac{dz}{2\pi\imath z^{m}}e^{(y+x_{l+1})z}\left((y+x_{l+1})z+z\partial_z+\nu+1\right)\frac{I_\nu(2\sqrt{yx_{l+1}}z)}{(\sqrt{yx_{l+1}}z)^\nu}\\
=&(m+\nu)\oint\frac{dz}{2\pi\imath z^{m}}e^{(y+x_{l+1})z}\frac{I_\nu(2\sqrt{yx_{l+1}}z)}{(\sqrt{yx_{l+1}}z)^\nu}.
\end{split}
\end{equation}
In the second line we pulled the differential operator into the integral since everything is absolutely integrable on the compact contour integral. In the third line we used the Bessel differential equation the Bessel function is solving and in the fourth line we used the fact that the function right of the bracket only depends on $\sqrt{yx_{l+1}}z$. In the last line we have integrated by parts. This result proofs the identity~\eqref{proof2.a}.

The counterpart of the calculation above for $\nu=-1/2$ is to check
\begin{equation}
[(\sqrt{z}\partial_z)^2 z^m]_{z=(\sqrt{y}\pm\sqrt{x_{l+1}})^2}=(\sqrt{y}\partial_y)^2(\sqrt{y}\pm\sqrt{x_{l+1}})^{2m}
\end{equation}
which is obviously true.

Collecting everything we have 
\begin{equation}\label{proof2.b}
\begin{split}
I_{ml}=&\int_0^\infty dy \int_{0}^\infty dr e^{-r}\omega(y)\left.\frac{(r-\partial_s{y}^{-\nu}\partial_{y}{y}^{\nu+1}\partial_{y})^{n-1}}{(n-1)!}\frac{1}{\mathcal{H}_\nu\omega(s)}\right|_{s=0}\\
&\times\left(\frac{\Gamma[\nu+1]}{\sqrt{\pi}\Gamma[\nu+1/2]}\int_{-1}^1\left(\prod_{k\neq m+1}\frac{x_k-y-x_{l+1}-2\sqrt{yx_{l+1}}t}{x_k-x_{m+1}}\right)(1-t^2)^{\nu-1/2}dt\right)
\end{split}
\end{equation}
for $\nu>-1/2$ and
\begin{equation}\label{proof2.b.2}
\begin{split}
I_{ml}=&\frac{1}{2}\int_0^\infty dy \int_{0}^\infty dr e^{-r}\omega(y)\left.\frac{(r-\partial_s(\sqrt{y}\partial_{y})^2)^{n-1}}{(n-1)!}\frac{1}{\mathcal{H}_{-1/2}\omega(s)}\right|_{s=0}\\
&\times\left(\prod_{k\neq m+1}\frac{x_k-y-x_{l+1}-2\sqrt{yx_{l+1}}}{x_k-x_{m+1}}+\prod_{k\neq m+1}\frac{x_k-y-x_{l+1}+2\sqrt{yx_{l+1}}}{x_k-x_{m+1}}\right)
\end{split}
\end{equation}
for $\nu=-1/2$.
We emphasize that the differential operator acts on everything what is right of it. Everything in the bracket is a polynomial of $y$ and can be expanded as $\sum_{j=0}^{n-1}c_j y^j$. Employing Theorem~\ref{thm:stat.Polya.b} and denoting the polynomials with
\begin{equation}\label{pol-proof.b}
\hat{p}_j(y)=\int_0^\infty dr e^{-r}\left.\frac{\Gamma(\nu+1)(r-\partial_t{y}^{-\nu}\partial_{y}{y}^{\nu+1}\partial_{y})^{n-1}}{(n-1)!j!\Gamma(\nu+j+1)}\frac{{y}^j}{\mathcal{H}_\nu\omega(t)}\right|_{t=0}
\end{equation}
we find
\begin{equation}
\begin{split}
I_{ml}=&\int_0^\infty dy \omega(y)\sum_{j=0}^{n-1}c_j \frac{j!\Gamma[\nu+j+1]}{\Gamma[\nu+1]}\hat{p}_j(y)
\end{split}
\end{equation}
which is true for all $\nu\in\{\pm1/2\}\cup\mathbb{N}_0$.
The bi-orthonormality of $\hat{p}_j(y)$ with $(\partial_{y}{y}^{\nu+1}\partial_{y}{y}^{-\nu})^j\omega(y)$ yields that only the coefficient $c_0$ survives. This coefficient is equal to $1$ only when $m=l$ and otherwise vanishes. This concludes the proof because for the kernel we only interchange the sum with the integrals.
\hfill$\square$
\vskip0.2cm

The polynomials in Eq.~\eqref{biorth.fixed.b} can be simplified again as well as the kernel when assuming the conditions of Corollary~\ref{cor:simpl.Polya.b} for the weight $\omega$.

\begin{corollary}[Simplification of  Theorem~\ref{thm:stat.fixed.b}]\label{cor:simpl.fixed.b}\

We assume the requirements of  Theorem~\ref{thm:stat.fixed.b} and the properties of $\omega$ in Corollary~\ref{cor:simpl.Polya.b}. Then the polynomials of the bi-orthonormal functions can be written as
 \begin{equation}\label{biorth.fixed.b.2}
 \begin{split}
 p_j(y')=&\oint\frac{dz'}{\pi\imath {z'}}\int_0^\infty dx\left(\frac{x}{y'z'}\right)^{\nu/2}\frac{I_\nu(2\sqrt{y'z'})K_\nu(2\sqrt{x})}{\mathcal{H}_\nu\omega(z')}\left(\prod_{l\neq j+1}\frac{x_l-x/z'}{x_l-x_{j+1}}\right)
 \end{split}
 \end{equation}
 while the weights are still given as in Eq.~\eqref{biorth.fixed.b}.
The kernel is either
 \begin{equation}\label{kernel.fixed.b.2}
 \begin{split}
 K_n(y',y)=&y^\nu\oint\frac{dz'}{\pi\imath {z'}}\int_0^\infty dx\left(\frac{x}{y'z'}\right)^{\nu/2}\frac{I_\nu(2\sqrt{y'z'})K_\nu(2\sqrt{x})}{\mathcal{H}_\nu\omega(z')}\\
 &\times\left(\sum_{j=1}^n\frac{\Gamma[\nu+1]}{\sqrt{\pi}\Gamma[\nu+1/2]}\int_{-1}^{1} \frac{\omega(y+x_{j}-2\sqrt{yx_{j}}t)}{(y+x_{j}-2\sqrt{yx_{j}}t)^{\nu}}(1-t^2)^{\nu-1/2}dt\prod_{l\neq j}\frac{x_l-x/z'}{x_l-x_j}\right)
 \end{split}
 \end{equation}
 for $\nu>-1/2$ and 
 \begin{equation}\label{kernel.fixed.b.3}
 \begin{split}
 K_n(y',y)=&\frac{1}{2\sqrt{y}}\oint\frac{dz'}{\pi\imath {z'}}\int_0^\infty dx\left(\frac{x}{y'z'}\right)^{\nu/2}\frac{I_\nu(2\sqrt{y'z'})K_\nu(2\sqrt{x})}{\mathcal{H}_\nu\omega(z')}\\
 &\times\left(\sum_{j=1}^n\left[ |\sqrt{y}-\sqrt{x_{j}}|\omega((\sqrt{y}-\sqrt{x_{j}})^2)+|\sqrt{y}+\sqrt{x_{j}}|\omega((\sqrt{y}+\sqrt{x_{j}})^2)\right]\prod_{l\neq j}\frac{x_l-x/z'}{x_l-x_j}\right)
 \end{split}
 \end{equation}
 for $\nu=-1/2$.
 The contour of $z'$ encircles the origin $z'=0$ counter clockwise.
\end{corollary}

{\bf Proof:}\\
This time we need the identity
\begin{equation}\label{proof.2.c}
\int_0^\infty dr e^{-r}\left.\frac{(r-\partial_t{y'}^{-\nu}\partial_{y'}{y'}^{\nu+1}\partial_{y'})^{n-1}}{(n-1)!}\frac{p(y')}{\mathcal{H}_\nu\omega(t)}\right|_{t=0}=\oint\frac{dz'}{\pi\imath {z'}}\int_0^\infty dx \left(\frac{x}{y'z'}\right)^{\nu/2}\frac{I_\nu(2\sqrt{y'z'})K_\nu(2\sqrt{x})}{\mathcal{H}_\nu\omega(z')}p\left(\frac{x}{z'}\right)
\end{equation}
for any polynomial $p(y')$ which has  maximally the order $n-1$ and any $\nu\in\{\pm1/2\}\cup\mathbb{N}_0$. The function $K_\nu$ is the modified Bessel function of the second kind. Again one can readily check this identity  by expanding the polynomial and evaluating the integrals.
\hfill$\square$
\vskip0.2cm

Once again we are sure that for particular cases the number of integrals can be reduced and even the sum can be rewritten into a contour integral, see Example~\ref{example.b.1}, as we have already seen for some P\'olya ensembles on $H_2$, cf.  Example~\ref{example.a.1}. However for a general P\'olya ensemble we were not able to simplify this result any further.

\begin{example}\label{example.b.1}\

We only want to consider the Laguerre ensemble with $\omega(x)=x^\nu e^{-x}$ which is already a new result. The weights in Eq.~\eqref{biorth.fixed.b} become essentially Bessel functions,
\begin{equation}\label{Lag.id.c}
 q_j(y)=\frac{\Gamma[\nu+1]}{\sqrt{\pi}\Gamma[\nu+1/2]}y^\nu\int_{-1}^{1} e^{-y-x_{j+1}+2\sqrt{yx_{j+1}}t}(1-t^2)^{\nu-1/2}dt=\Gamma[\nu+1]\left(\frac{y}{x_{j+1}}\right)^{\nu/2} I_\nu(2\sqrt{yx_{j+1}})e^{-y-x_{j+1}}
\end{equation}
with $\nu>-1/2$. For $\nu=-1/2$ the right hand side still holds as can be easily checked with the identification $\Gamma[1/2]I_{-1/2}(2\sqrt{x})=x^{-1/4}\cosh(2\sqrt{x})$.

The weights $q_j$ are entire in $x_{j+1}$. Hence we can rewrite the sum in the kernel~\eqref{kernel.fixed.b.2} as a contour integral for $\nu>-1/2$
 \begin{equation}
 \begin{split}
&\sum_{j=1}^n\frac{\Gamma[\nu+1]}{\sqrt{\pi}\Gamma[\nu+1/2]}y^\nu\int_{-1}^{1} \frac{\omega(y+x_{j}-2\sqrt{yx_{j}}t)}{(y+x_{j}-2\sqrt{yx_{j}}t)^{\nu}}(1-t^2)^{\nu-1/2}dt\prod_{l\neq j}\frac{x_l-x/z'}{x_l-x_j}\\
=&-\Gamma[\nu+1]\oint\frac{dz}{2\pi\imath}\frac{1}{z-x/z'}\left(\frac{y}{z}\right)^{\nu/2} I_\nu(2\sqrt{yz})e^{-y-z}\prod_{l=1}^n\frac{x_l-x/z'}{x_l-z}
 \end{split}
 \end{equation}
 where $z$ only encircles the poles $x_1,\ldots,x_n$ (but not $x/z'$) counter clockwise, and similar for $\nu=-1/2$ which does not change the result.
For calculating the polynomials we need the identity
\begin{equation}\label{Lag.id.b}
\begin{split}
\oint\frac{dz'}{\pi\imath {z'}}\int_0^\infty dx\left(\frac{x}{y'z'}\right)^{\nu/2}I_\nu(2\sqrt{y'z'})K_\nu(2\sqrt{x})e^{z'}\left(\frac{x}{z'}\right)^m=&\sum_{j=0}^m\frac{m!\Gamma[\nu+m+1]}{j!(m-j)!\Gamma[\nu+j+1]}{y'}^j\\
=&m!\oint \frac{dz}{2\pi\imath z^{\nu+1}(1-z)^{m+1}}\exp[-y'(1-z^{-1})]\\
=&\int_0^\infty dx'\left(\frac{x'}{y'}\right)^{\nu/2}e^{-x'-y'}I_\nu(2\sqrt{x'y'}){x'}^m
\end{split}
\end{equation}
for any $m=0,\ldots,n-1$ which is true for any $\nu\in\mathbb{N}_0$. In the second line we close the contour around origin counter clockwise but do not enclose the pole at $z=1$. We can extend the identity~\eqref{Lag.id.b} to $\nu=\pm1/2$ by explicit evaluation of the integral over $x'$. Hence the polynomials~\eqref{biorth.fixed.b.2} have the simple form
 \begin{equation}
 \begin{split}
 p_j(y')=&\int_0^\infty dx'\left(\frac{x'}{y'}\right)^{\nu/2}e^{-x'-y'}I_\nu(2\sqrt{x'y'})\left(\prod_{l\neq j+1}\frac{x_l-x'}{x_l-x_{j+1}}\right)
 \end{split}
 \end{equation}
 and the kernel~\eqref{kernel.fixed.b.2} becomes then
 \begin{equation}
 \begin{split}
 K_n(y',y)=&-\Gamma[\nu+1]\int_0^\infty dx'\oint\frac{dz}{2\pi\imath}\frac{1}{z-x'}\left(\frac{yx'}{y'z}\right)^{\nu/2} I_\nu(2\sqrt{x'y'})I_\nu(2\sqrt{yz})e^{-x'-y'-y-z}\prod_{l=1}^n\frac{x_l-x'}{x_l-z}
 \end{split}
 \end{equation}
 which is true for all $\nu\in\{\pm1/2\}\cup\mathbb{N}_0$. The contour of $z$ only encircles the points $x_1,\ldots,x_n$ (but not $x'$) counter clockwise. The result for the kernel is completely new while there was already a formula derived for the polynomials in terms of an integral over a hypergeometric function in~\cite{Forrester:2013}. 
\end{example}

\subsection{Statistics with Polynomial Ensembles on $M=H_1,H_4,M_\nu$}\label{sec:evstat.poly.b}

Next we want to consider the case of a convolution of a P\'olya ensemble on $M=H_1,H_4,M_\nu$ and a general polynomial ensemble on $M$.

\begin{theorem}[Eigenvalue Statistics of P\'olya Ensembles with Polynomial Ensembles]\label{thm:stat.poly.b}\

 Consider the joint probability density~\eqref{jpdf.poly.b}  where the polynomial ensemble of $X_2\in M=H_1,H_4,M_\nu$ can be described by the bi-orthonormal functions $\{\tilde{p}_j,w_{j+1}\}_{j=0,\ldots,n-1}$ and its kernel is $\tilde{K}_n(y',y)=\sum_{j=0}^{n-1}\tilde{p}_j(y')w_{j+1}(y)$. The pair of bi-orthonormal functions corresponding to $Y=X_1+X_2$ is
 \begin{equation}\label{biorth.poly.b}
 \{p_j,q_j\}_{j=0,\ldots,n-1}=\left\{\int_{0}^\infty dr e^{-r}\left.\frac{(r-\partial_t{y'}^{-\nu}\partial_{y'}{y'}^{\nu+1}\partial_{y'})^{n-1}}{(n-1)!}\frac{\tilde{p}_j(y')}{\mathcal{H}_\nu\omega(t)}\right|_{t=0}\,,\,\omega\ast_\nu w_{j+1}(y)\right\}_{j=0,\ldots,n-1}
 \end{equation}
and the corresponding kernel is
 \begin{equation}\label{kernel.poly.b}
 \begin{split}
 K_n(y',y)=&y^\nu\int_{0}^\infty dx\int_{0}^\infty dr e^{-r}\left(\frac{\Gamma[\nu+1]}{\sqrt{\pi}\Gamma[\nu+1/2]}\int_{-1}^1\frac{\omega(x+y-2\sqrt{xy} t)}{(x+y-2\sqrt{xy} t)^{\nu}}(1-t^2)^{\nu-1/2}dt\right)\\
 &\times\left.\frac{(r-\partial_t{y'}^{-\nu}\partial_{y'}{y'}^{\nu+1}\partial_{y'})^{n-1}}{(n-1)!}\frac{\tilde{K}_n(y',x)}{\mathcal{H}_\nu\omega(t)}\right|_{t=0}
 \end{split}
 \end{equation}
 for $\nu>-1/2$ and
 \begin{equation}\label{kernel.poly.b.b}
 \begin{split}
 K_n(y',y)=&\frac{1}{2\sqrt{y}}\int_{0}^\infty dx\int_{0}^\infty dr e^{-r}\left[|\sqrt{x}-\sqrt{y}|\omega((\sqrt{x}-\sqrt{y})^2)+|\sqrt{x}+\sqrt{y}|\omega((\sqrt{x}+\sqrt{y})^2)\right]\\
 &\times\left.\frac{(r-\partial_t{y'}^{-\nu}\partial_{y'}{y'}^{\nu+1}\partial_{y'})^{n-1}}{(n-1)!}\frac{\tilde{K}_n(y',x)}{\mathcal{H}_\nu\omega(t)}\right|_{t=0}
 \end{split}
 \end{equation}
 for $\nu=-1/2$.
\end{theorem}

{\bf Proof:}\\
The proof works along the same line as for the proof of Theorem~\ref{thm:stat.fixed.b} and the difference to Eq.~\eqref{proof2.b} is that we have to integrate over $x$ and replace the product by the polynomial $\tilde{p}_m$, too, i.e.
\begin{equation}
\begin{split}
I_{ml}=&\int_0^\infty dx\int_0^\infty dy \int_{0}^\infty dr e^{-r}\omega(y) w_{l+1}(x)\left.\frac{(r-\partial_s{y}^{-\nu}\partial_{y}{y}^{\nu+1}\partial_{y})^{n-1}}{(n-1)!}\frac{1}{\mathcal{H}_\nu\omega(s)}\right|_{s=0}\\
&\times\left\{\begin{array}{cl} \displaystyle\frac{\Gamma[\nu+1]}{\sqrt{\pi}\Gamma[\nu+1/2]}\int_{-1}^1\tilde{p}_m(x+y+2\sqrt{xy}t)(1-t^2)^{\nu-1/2}dt, & {\rm for}\ \nu>-1/2, \\ \displaystyle\frac{1}{2}\left[\tilde{p}_m((\sqrt{x}-\sqrt{y})^2)+\tilde{p}_m((\sqrt{x}+\sqrt{y})^2)\right], & {\rm for}\ \nu=-1/2. \end{array}\right.
\end{split}
\end{equation}
We again expand what is in the bracket which yields the sum $\sum_{j=0}^md_j(x)y^j$ with $d_j(x)$ polynomials of $x$. The integral over $r$ with the differential operator generates the polynomials $\hat{p}_j(y)$ from the monomial $y^j$, see Eq.~\eqref{pol-proof.b}. These polynomials are bi-orthonormal to $(\partial_{y}{y}^{\nu+1}\partial_{y}{y}^{-\nu})^j\omega(y)$ such that the integral over $y$ selects only the coefficient $d_0(x)=\tilde{p}_m(x)$. Thus the integral becomes
\begin{equation}
I_{ml}=\int_0^\infty dx w_{l+1}(x)\tilde{p}_m(x)=\delta_{ml}
\end{equation}
which proves the bi-orthonormality of $\{p_j,q_j\}_{j=0,\ldots,n-1}$. The kernel is again a direct consequence since we have only to interchange the finite sum with the integrals.
\hfill$\square$
\vskip0.2cm

Again we conclude our general theorem with a simplification when assuming the conditions of $\omega$ as in Corollary~\ref{cor:simpl.Polya.b}.

\begin{corollary}[Simplification of  Theorem~\ref{thm:stat.fixed.b}]\label{cor:simpl.poly.b}\

We assume the requirements of  Theorem~\ref{thm:stat.poly.b} and the properties of $\omega$ in Corollary~\ref{cor:simpl.Polya.b}. Then the bi-orthonormal functions can be written as
 \begin{equation}\label{biorth.poly.b.2}
 \begin{split}
 \{p_j,q_j\}_{j=0,\ldots,n-1}=&\biggl\{\oint\frac{dz'}{2\pi\imath {z'}}\int_0^\infty dx\left(\frac{x}{y'z'}\right)^{\nu/2}\frac{I_\nu(2\sqrt{y'z'})K_\nu(2\sqrt{x})}{\mathcal{H}_\nu\omega(z')}\tilde{p}_j\left(\frac{x}{z'}\right)\,,\,\omega\ast_\nu w_{j+1}(y)\biggl\}_{j=0,\ldots,n-1}.
 \end{split}
 \end{equation}
and the kernel as
 \begin{equation}\label{kernel.poly.b.2}
 \begin{split}
 K_n(y',y)=&y^\nu\int_{0}^\infty d\hat{y}\oint\frac{dz'}{2\pi\imath {z'}}\int_0^\infty dx\left(\frac{x}{y'z'}\right)^{\nu/2}\frac{I_\nu(2\sqrt{y'z'})K_\nu(2\sqrt{x})}{\mathcal{H}_\nu\omega(z')}\\
 &\times\left(\frac{\Gamma[\nu+1]}{\sqrt{\pi}\Gamma[\nu+1/2]}\int_{-1}^1\frac{\omega(\hat{y}+y-2\sqrt{\hat{y}y} t)}{(\hat{y}+y-2\sqrt{\hat{y}y} t)^{\nu}}(1-t^2)^{\nu-1/2}dt\right)\tilde{K}_n\left(\frac{x}{z'},\hat{y}\right)
 \end{split}
 \end{equation}
 for $\nu>-1/2$ and
 \begin{equation}\label{kernel.poly.b.3}
 \begin{split}
 K_n(y',y)=&\frac{1}{2\sqrt{y}}\int_{0}^\infty d\hat{y}\oint\frac{dz'}{2\pi\imath {z'}}\int_0^\infty dx\left(\frac{x}{y'z'}\right)^{\nu/2}\frac{I_\nu(2\sqrt{y'z'})K_\nu(2\sqrt{x})}{\mathcal{H}_\nu\omega(z')}\\
 &\times\left[ |\sqrt{y}-\sqrt{\hat{y}}|\omega((\sqrt{y}-\sqrt{\hat{y}})^2)+|\sqrt{y}+\sqrt{\hat{y}}|\omega((\sqrt{y}+\sqrt{\hat{y}})^2)\right]\tilde{K}_n\left(\frac{x}{z'},\hat{y}\right)
 \end{split}
 \end{equation}
 for $\nu=-1/2$.
 The contour of $z'$ encircles the origin $z'=0$ counter clockwise.
\end{corollary}

{\bf Proof:}\\
The identity~\eqref{proof.2.c} yields these statements.
\hfill$\square$
\vskip0.2cm

As we will show in the next example also here the number integrals can be still reduced for specific ensembles.

\begin{example}\label{example.b.2}\

Again we consider the Laguerre ensemble with $\omega(x)=x^\nu e^{-x}$. We make use of the results~\eqref{Lag.id.c} and ~\eqref{Lag.id.b}  for a shift with fixed matrix in $M$ since they still apply in a slightly modified way. Thus the bi-orthonormal functions are
\begin{equation}
\begin{split}
 p_j(y')=&\int_0^\infty dx'\left(\frac{x'}{y'}\right)^{\nu/2}e^{-x'-y'}I_\nu(2\sqrt{x'y'})\tilde{p}_j(x'),\\
 q_j(y)=&\Gamma[\nu+1]\int_0^\infty d\hat{y}\left(\frac{y}{\hat{y}}\right)^{\nu/2}  I_\nu(2\sqrt{y\hat{y}})e^{-y-\hat{y}}\tilde{q}_j(\hat{y}).
\end{split}
\end{equation}
The kernels~\eqref{kernel.poly.b.2} and~\eqref{kernel.poly.b.3}  simplify then to the unified expression
 \begin{equation}
 \begin{split}
 K_n(y',y)=&\Gamma[\nu+1]\int_0^\infty dx'\int_0^\infty d\hat{y} \left(\frac{yx'}{y'\hat{y}}\right)^{\nu/2}e^{-x'-y'-y-\hat{y}}I_\nu(2\sqrt{x'y'}) I_\nu(2\sqrt{y\hat{y}})\tilde{K}_n\left(x',\hat{y}\right).
 \end{split}
 \end{equation}
 This is also a completely new result.
\end{example}

\section{Conclusions}\label{sec:conclusio}

We derived closed expressions for the bi-orthogonal functions and kernels for the eigenvalues/squared singular values, respectively, of three situations of random matrices. The first case is for a general P\'olya ensemble without a shift on either the Hermitian antisymmetric matrices $H_1$, the Hermitian matrices $H_2$, the Hermitian anti-self-dual matrices $H_4$ or the complex rectangular matrices $M_\nu$.  The other two cases we considered are the eigenvalue/squared singular value statistics of the P\'olya ensemble added by a either a fixed matrix or a random matrix drawn from a polynomial ensemble on the same space as the P\'olya ensemble, see~\cite{Kuijlaars:2014,Kuijlaars:2014b,Kieburg:2015,Liu:2016,Claeys:2016,Liu:2017,Akemann:2017}. All results hold for finite matrix dimension. The next step would be to analyse the asymptotic limits of these statistics in the limit of large matrices. We are sure that the expressions in terms of a very small number of integrals to be performed provide a good starting point to study a broad class of ensembles. Regarding this point we want to underline that many classical random matrix ensembles fall into one of the three kinds of P\'olya ensembles which were discovered~\cite{Kieburg:2016a,Kieburg:2016b,Kuijlaars:2016a,Kieburg:2017}. However the range of the P\'olya ensembles goes far beyond these classical examples and even yield highly non-trivial cases~\cite{Kieburg:2016a,Kieburg:2016b,Kieburg:2017}. Hereby we also want to emphasize that the class of P\'olya ensembles usually yield random matrix ensembles which have not the simple potential form $P(X)\propto\exp[-\tr V(X)]$ as studied in the broad literature, e.g. see~\cite{Borodin:1999}. Thus many new things and phenomena may occur.

The technique we have used is the method of spherical functions and transforms from harmonic analysis on Lie groups, see~\cite{Helgason:2000}. It proved again as a very effective tool to deal with convolutions regardless whether they are of multiplicative nature as in~\cite{Kieburg:2016a,Kieburg:2016b} or of additive one as in~\cite{Kuijlaars:2016a,Kieburg:2017} or here. This tool seem to be also suitable studying stochastic processes on matrix spaces. Thus another direction of further investigation could be the analysis of stochastic processes at finite matrix dimension and for a finite number matrices involved and their various limits where the rate of convergence might be tractable.

\section*{Acknowledgements}

I am grateful for the intensive discussions with Holger K\"osters. Special thanks also goes to  Tomasz Checinski and Peter Forrester who gave me detailed feedback on a pre-version of the present work. Moreover I want to thank  Gernot Akemann, Tom Claeys, Arno Kuijlaars, Jesper Ipsen, and Eugene Strahov
for helpful comments on this topic. I acknowledge support by the grant AK35/2-1 ``Products of Random Matrices" of the German research council (DFG).

\end{document}